\documentclass[a4paper,11pt,abstracton]{scrartcl}
\pdfoutput=1
\usepackage{multicol}
\usepackage{color}
\usepackage{xspace}

\usepackage{amsfonts}
\usepackage{amsthm}
\usepackage{amsmath}
\usepackage{amssymb}
\usepackage{graphicx}%

\usepackage{graphicx}
\usepackage{float}
\usepackage[autolanguage]{numprint}  
\usepackage{hyperref} 
 \usepackage{booktabs} 
\usepackage{fixmath,amsmath,amssymb,amsfonts} 

\usepackage{subfig}
\usepackage{float}

\newcommand{\eff}{{{\text{eff}}}}

\newcommand {\pnull}{\ensuremath{{P}^D}}

\newcommand{\ve}{\ensuremath{{\v}^\epsilon}}
\newcommand{\vecomp}[1]{\ensuremath{{v}_{#1}^\epsilon}}
\newcommand{\pe}{\ensuremath{p^\epsilon}}

\newcommand{\cbl}{\ensuremath{C^{2,bl}_1}}
\newcommand{\cblpi}{\ensuremath{C^{bl}_\pi}}

\def\ve{\mathbf{v}^\ep}
\def\ep{\varepsilon}

\def\p{\partial}

\def\u2{{u^\ep \over \ep^2 }}
\def\u3{{\displaystyle {\bar u}^\ep \over \ep^2 }}
\def\p{\partial}
\def\eop{~\vrule width 0.4 pc height 0.6 pc depth 0.1 pc}

\def\ep{\varepsilon}

\def\inf{\infty}

\def\O{\Omega}
\def\Oe{\Omega^\ep}

\def\u2{{u^\ep \over \ep^2 }}
\def\u3{{{\bar u}^\ep \over \ep^2 }}

\def\z1{\zeta^\ep_1}
\def\C1{C_{i,\zeta }}
\def\C2{C_{1,\zeta }}

\def\p{\partial}

\newtheorem{theorem}{Theorem}

\newtheorem{corollary}[theorem]{Corollary}

\newtheorem{proposition}[theorem]{Proposition}
\newtheorem{remark}[theorem]{Remark}

\numberwithin{equation}{section}

\title{Effective  interface conditions  for the forced infiltration of a viscous fluid into a porous medium using homogenization}

\author{T. Carraro\thanks{\texttt{thomas.carraro@iwr.uni-heidelberg.de}, Institute for Applied Mathematics, Heidelberg University, 69120
Heidelberg, Germany. TC was supported by the German Research Council (DFG) through project
 ``Modellierung, Simulation und Optimierung der Mikrostruktur mischleitender
 SOFC-Kathoden'' (RA 306/17-2).} 
 \and C. Goll\thanks{\texttt{christian.goll@iwr.uni-heidelberg.de}, Institute for Applied Mathematics, Heidelberg University, 69120
Heidelberg, Germany.} 
 \and A. Marciniak-Czochra\thanks{\texttt{anna.marciniak@iwr.uni-heidelberg.de}, Bioquant, Heidelberg University, 69120 Heidelberg, Germany. AM-C  was supported by ERC Starting Grant "Biostruct" No. 210680 and Emmy Noether Programme of German Research Council (DFG).} 
 \and A. Mikeli\'c\thanks{\texttt{andro.mikelic@univ-lyon1.fr} \textit{Corresponding author} Universit\'e de Lyon, Lyon, F-69003, France; Universit\'e
 Lyon 1, Institut Camille Jordan, UMR 5208, B\^at. Braconnier,  43, Bd du 11
 novembre 1918, 69622 Villeurbanne Cedex,
 France. The research of A.M. was partially supported by the  Programme Inter Carnot Fraunhofer from BMBF (Grant 01SF0804) and ANR. Research visits of A.M. to the Heidelberg University were supported in part by the Romberg professorship at IWR, Heidelberg University, 2011--1013.}}%

\begin{document}

 \maketitle
\begin{abstract}  It is generally accepted that the effective velocity of a viscous flow over
a porous bed satisfies the Beavers-Joseph slip law. To the contrary, in the case of a forced infiltration of a viscous fluid into a porous medium the interface law has been a subject of controversy.
 In this paper, we prove rigorously that the effective { interface } conditions are: (i) the continuity of the normal effective velocities; (ii) zero Darcy's pressure and (iii) a { given } slip velocity.
 { The effective tangential slip velocity is  calculated from the boundary layer and depends only on the pore geometry. In the next order of approximation, we derive a pressure slip law}.  An independent confirmation of the analytical results using direct numerical simulation of the flow at the microscopic level is given, as well.
\end{abstract}
\textsf{\textbf{\small Keywords}} Interface conditions, pore scale simulation, pressure slip law, slip velocity, boundary layers, homogenization
%
%

\section{Introduction}
The purpose of this paper is to derive rigorously the interface conditions governing the infiltration of a viscous fluid into a porous medium.

We start from an incompressible 2D flow of a Newtonian fluid penetrating a porous medium. { At the pore scale, the}  flow is described by the stationary Stokes system, both, in the unconstrained fluid part and in the pore space. The upscaling of the Stokes system in a porous medium yields Darcy's law as the effective momentum equation, valid at every point of the porous medium. The Stokes system and the Darcy equation are very different PDEs and need to be coupled at the interface of the fluid and the porous medium. The resulting system should be an approximation of the starting first principles with  error estimate in the term of the dimensionless pore size $\ep$, being the ratio of the characteristic pore size and the macroscopic domain length.

There is vast literature on modeling interface conditions between a free flow and a porous medium. Most of the references focus on flows which are tangential to the porous medium. In such situation, the free fluid velocity is much larger than the Darcy velocity in the porous medium. The corresponding interface condition is the slip law by Beavers and Joseph. It was deduced from the experiment in \cite{BJ}, then discussed and simplified into a generally used form in \cite{SAF} and justified through numerical simulations  of pore level Navier-Stokes equations in  \cite{SahKav92}, \cite{K95} and \cite{CGMCM:13}. A rigorous justification of the slip law by Beavers and Joseph, starting from the pore level first principles, was provided by J\"ager and Mikeli\'c  in \cite{JM00}, using a combination of homogenization and boundary layers techniques. The slip law is supplemented by the pressure jump law, what was noticed in \cite{JMN01} and rigorously derived in \cite{AMCAM2011}. A corresponding numerical validation
by solving the Stokes equation at the pore scale has been recently presented in \cite{CGMCM:13}.

Infiltration into a porous medium corresponds to a different situation, because in this case the free fluid velocity and the Darcy velocity are of the same order. We refer to the article by Levy and Sanchez-Palencia \cite{LSP75}. They classify the physical situation as "Case B: The pressure gradient on the side of the porous body at the interface is normal to it". { In the "Case B" the pressure gradient in the porous medium is much larger than in the free fluid.} Using
 { an order-of-magnitude analysis}, in  \cite{LSP75} it was concluded that the effective interface conditions have to satisfy
\begin{gather} \mathbf{u}^{eff} \cdot \mathbf{n} = \mathbf{u}^{D} \cdot \mathbf{n} \quad \mbox{ and } \quad P^D =\mbox{ a constant},
\label{Interfc1}
\end{gather}
where $\{ \mathbf{u}^{D} , P^D \}$ are the Darcy velocity and the pressure and $\mathbf{u}^{eff}$ is the unconfined fluid velocity. { Note that the interface conditions (\ref{Interfc1}) were obtained for low Reynolds number flows.}

In order to couple the Stokes system in the free fluid domain with the Darcy equation in the porous medium, the conditions (\ref{Interfc1}) are not sufficient. One more condition is needed. In \cite{LSP75} an intermediate boundary layer was introduced and existence of an effective slip velocity at the interface was postulated. However, the article \cite{LSP75} did not provide the slip velocity. It was limited to a model of macroscopic isotropy, where the slip is equal to zero. Therefore, zero tangential velocity is assumed.

A rigorous mathematical study of the interface conditions between a free fluid and a porous medium was initiated in \cite{JaMi2}. Our analysis  repose on the boundary layers constructed there. For reviews of the models and techniques we refer to \cite{QuartRev09} and \cite{JaegMik09}.

We note that in a number of  articles devoted to numerical simulations, the porous part was modeled using the Brinkman-extended Darcy law. We refer to \cite{Quat}, \cite{HWNW06}, \cite{IL}, \cite{NHWRW05}, \cite{Yuetc07} and references therein. In such setting, the authors used general interface conditions introduced by Ochoa-Tapia and Whitaker in \cite{OTW1:95}. They consists of (i) continuity of the velocity and (ii) complex jump relations for the stresses, containing several parameters to be fitted. We recall that the viscosity in the Brinkman equation is not known and the use of it seems to be justified only in the case of a high porosity (see the discussion in \cite{Nield09}). Furthermore, Larson and Higdon undertook a detailed numerical simulation of two configurations (axial and transverse) of a shear flow over a porous medium in \cite{LH86}.  Their conclusion was that a macroscopic model based on Brinkman's equation gives ``reasonable predictions for the rate of decay of the mean
velocity for certain simple geometries, but fails for to predict the correct behavior for media anisotropic in the plane normal to the flow direction''.  An approach using the thermodynamically constrained averaging theory was presented in \cite{JRHGM12}. Darcy-Navier-Stokes coupling yields also an increasing interest from the side of numerical analysis, see  \cite{Yot}, \cite{Riviere05} and \cite{QuartRev09} and references therein.

In this paper, we provide a rigorous derivation of the filtration equation and the interface condition explained above from the pore scale level description based on first principles. Our derivation follows the general homogenization and boundary layers approach from \cite{JaMi2}.  The necessary results on boundary layers and very weak solutions to the Stokes system will be recalled in the proofs of the main results.

{ In our work we have used a finite element method to obtain a numerical confirmation of the analytical results. The numerical study of the convergence rates of the macroscopic problems and effective interface conditions is a difficult task for the reasons explained in the following. The first difficulty is the numerical solution of the microscopic problem used as reference. The geometry of the porous part has to be resolved with high accuracy. In addition, the microscopic solution in the vicinity of the surface of the porous medium has large gradients, that can only be approximated by a boundary layer as shown in this work. The accuracy needed by the resolution of the interface and porous part requires high performance computing. We have reduced the computational costs by considering for our test cases a problem with periodic geometry and periodic boundary conditions. We could thus reduce our computations to one column of inclusions in the porous part. Nevertheless, even in our simplified
example problem all the computations must be performed with high accuracy. The reason is that the homogenization errors, especially in the estimates based on correction terms, are small in comparison with numerical errors even for simulations with millions of degrees of freedom. A further difficulty is that to numerically check the estimates we have to solve several auxiliary problems that are coupled. Therefore the numerical precision of one problem influences the precision of the other ones. Due to the complexity of the microscopic flow and the boundary layers, strategies for local mesh adaptivity to reduce the computations of the norms in the estimates are not effective. We could nevertheless apply a goal oriented adaptive method, based on the dual weighted residual (DWR) method \cite{BeckeR:2001}, to calculate some constants needed for the estimates, increasing the overall accuracy of our numerical tests.}

 The paper is organized as follows: In Section \ref{StatPb} we formulate the starting microscopic problem and the resulting effective equations. We provide theorems on error estimates of the model approximation. In Section \ref{Numerics} we give a numerical confirmation of the analytical results based on finite element computations. Sections \ref{Proofs1}--\ref{Proofs3} contain the corresponding proofs.


\section{Problem setting and main results}\label{StatPb}
\subsection{Definition of the geometry} \label{ExpDarcy}


Let $L, h$ and $H$ be positive real numbers. We consider a two dimensional periodic porous medium $\Omega_2 =
(0,L)\times (-H, 0)$ with a periodic arrangement of the pores. The formal
description goes along the following lines: \vskip0pt
First, we define the geometrical
structure inside the unit cell $Y = (0,1)^2$. Let
$Y_s$ (the solid part) be a closed strictly included subset of $\bar{Y}$, and $Y_F =
Y\backslash Y_s$ (the fluid part). Now we make a periodic
repetition of $Y_s$ all over $\mathbb{R}^2$ and set $Y^k_s = Y_s + k $, $k
\in \mathbb{Z}^2$. Obviously, the resulting set $E_s = \bigcup_{k \in
\mathbb{Z}^2} Y^k_s$ is a closed subset of $ \mathbb{R}^2$ and $E_F =  \mathbb{R}^2
\backslash E_s$ in an open set in $ \mathbb{R}^2$.  We suppose that
 $Y_s$ has a boundary of class $ C^\infty$, which is locally located on
one side of their boundary. Obviously,  $ E_F $ is connected and
$E_s$ is not.  \vskip3pt Now we notice that $\Omega_2$ is covered with a
regular mesh of size $ \varepsilon$, each cell being a cube
$Y^{\varepsilon}_i$, with $1 \leq i \leq N(\varepsilon) = \vert
\Omega_2 \vert  \varepsilon^{-2} [1+ o(1)]$. Each cube
$Y^{\varepsilon}_i$ is homeomorphic to $Y$, by linear homeomorphism
$\Pi^{\varepsilon}_i$, being composed of translation and a homothety
of ratio $1/ \varepsilon$.

We define
$\displaystyle
Y^{\varepsilon}_{S_i} = (\Pi^{\varepsilon}_i)^{-1}(Y_s)\qquad
\hbox{ and }\quad Y^{\varepsilon}_{F_i} =
(\Pi^{\varepsilon}_i)^{-1}(Y_F).
$
For sufficiently small $\varepsilon > 0 $ we consider the set
$ \displaystyle T_{\varepsilon} = \{k \in  \mathbb{Z}^2  \vert  Y^{\varepsilon}_{S_k} \subset
\Omega_2 \} $
and define
$$
O_{\varepsilon} = \bigcup_{k \in T_{\varepsilon}}
Y^{\varepsilon}_{S_k} , \quad S^{\varepsilon} = \partial
O_{\varepsilon}, \quad \Omega^{\varepsilon}_2 = \Omega_2 \backslash
O_{\varepsilon} =\O_2 \cap \ep E_F$$ Obviously, $\partial
\Omega^{\varepsilon}_2 = \partial \Omega_2 \cup S^{\varepsilon}$. The
domains $O_{\varepsilon}$ and $\Omega^{\varepsilon}_2 $ represent,
respectively, the solid and fluid parts of the porous medium
$\Omega$. For simplicity, we suppose $L/\varepsilon , H/\varepsilon , h/\varepsilon \in \mathbb{N}$.

We set $\Sigma =(0,L) \times \{ 0\} $, $\Omega_1 = (0,L)\times (0,h)$ and $\O = (0,L) \times (-H,h)$.
Furthermore, let $\O^\ep = \Oe_2 \cup \Sigma \cup \O_1 $ and $\O = \O_2 \cup \Sigma \cup \O_1 $.

\begin{figure}
\center
 \subfloat[\label{afoto}]{
\resizebox{0.4\textwidth}{!}{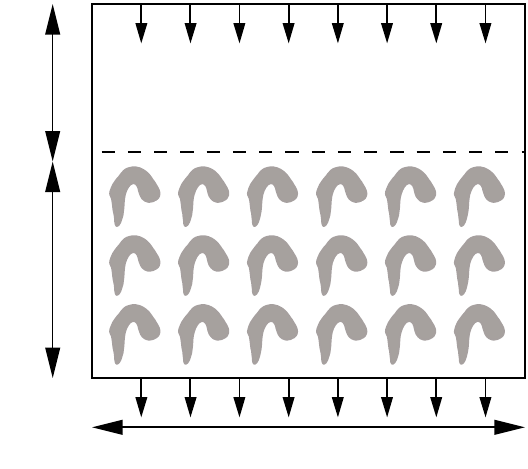}}
\hspace{1cm}
\subfloat[\label{subfig.unit cell}]{
\raisebox{1.35cm}{\resizebox{0.25\textwidth}{!}{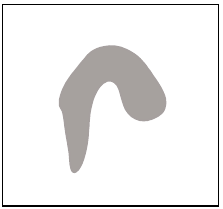}}}
\caption{Sketch of the geometry
\protect\subref{afoto} the periodicity cell $Y$~ \protect\subref{subfig.unit cell}.}
\end{figure}


\subsection{The microscopic model}\label{MicroP}

Having
defined the geometrical structure of the porous medium, we precise
the flow problem.

We consider the  slow viscous incompressible
flow of a single fluid through a porous medium. The flow is caused by the fluid injection at the boundary $\{ x_2 =h \}.$ We suppose the
no-slip condition at the boundaries of the pores (i.e. the filtration through a rigid porous medium). Then, the flow is described by the following non-dimensional steady Stokes system in $\Omega^\varepsilon $
(the fluid part of the porous medium $\Omega$):
\begin{subequations}\label{equ.micro}
 \begin{gather} -
 \Delta \mathbf{v}^\ep + \nabla p^{\ep} = 0 \qquad \hbox{ in } \quad \Oe ,
\label{1.3} \\ \; \mbox{ div }  \, \ve = 0 \qquad \hbox{ in } \quad \Oe , \qquad \int_{\O_1} p^\ep \ dx =0,
\label{1.4} \\ \ve  =0 \quad \hbox{on } \quad \p \Oe \setminus \O 
, \qquad \{ \ve , p^\ep \} \quad \hbox{ is }
L-\hbox{periodic in } \; x_1 ,  \label{1.5} \\
\mathbf{v}^\ep |_{x_2 =h} = \mathbf{v}^D , \quad 
{v}^\ep_2 |_{x_2 =-H} = g, \quad 
\frac{\p v^\ep_1}{\p x_2} |_{x_2 =-H} = 0 .\label{1.8}
\end{gather}
\end{subequations}
Such flow is possible only under  the following compatibility condition
\begin{equation}\label{Compa}
   L U_B = \int^L_0 g(x_1 ) \ dx_1 = \int^L_0 v^D_2 (x_1 ) \ dx_1 .
\end{equation}
With the assumption on the geometry from section \ref{ExpDarcy}, condition (\ref{Compa}) and for $\mathbf{f}\in C^\infty ({\overline \Omega})^2$, $\mathbf{v}^D \in C^\infty [0,L]^2$ and $g\in C^\infty [0,L]$, problem (\ref{1.3})-(\ref{1.8}) admits a unique solution $\{ \mathbf{v}^\ep , p ^\ep \} \in C^\infty ({\overline \Oe})^3$, for all $\ep >0$.

Our goal is to study behavior of solutions to system (\ref{1.3})-(\ref{1.8}), when $\ep \to 0$. In such limit the equations in $\Omega_1$ remain unchanged and in $\Oe_2$  the Stokes system is upscaled to Darcy's equation posed in $\Omega_2$. Our contribution is the derivation of the interface condition, linking these two systems.
\subsection{The boundary layers and effective coefficients}
 \vskip0pt
Let the permeability tensor $K$ be given by
\begin{equation}\label{1.58}
   K_{ij} = \int_{Y_F} \nabla_y \mathbf{w}^i : \nabla_y \mathbf{w}^j \ dy =\int_{Y_F} w^i_j \ dy ,
\; 1\leq i,j \leq 2.
\end{equation}
where
$\hbox{ for } 1\leq i\leq 2 ,\,  \ \{ \mathbf{w}^i ,
\pi^i \} \in H^1_{per} (Y_F)^2 \times L^2 (Y_F), \,  \int_{Y_F}  \pi^i (y ) \, dy =0, \hbox{ are solutions to
} $
\begin{equation}\label{1.57Cell}
   \left\{ \begin{matrix}
   \hfill - \Delta_y \mathbf{w}^i (y) + \nabla_y \pi^i (y ) = \mathbf{e}^i &
\hbox{ in } \ Y_F \cr \hfill \hbox{ div}_y \mathbf{w}^i (y) =0  & \hbox{ in
} \ Y_F \cr \hfill \mathbf{w}^i (y ) =0 & \hbox{ on } \, (\partial Y_F
\setminus
\partial Y). \cr
 \end{matrix}  \right.
\end{equation}
 Obviously, these problems always admit  unique
solutions and $K$ is a symmetric positive definite matrix (the dimensionless permeability tensor).
In addition
\begin{equation}\label{Ksurf}
  K_{j2}=K_{2j}= \int^1_0 w^j_2 (y_1 , 0) \ dy_1 .
\end{equation}
\vskip3pt
In order to formulate the result we need the  viscous boundary layer problem connecting free fluid flow and a porous medium flow:
\begin{figure}
\center
\resizebox{0.25\textwidth}{!}{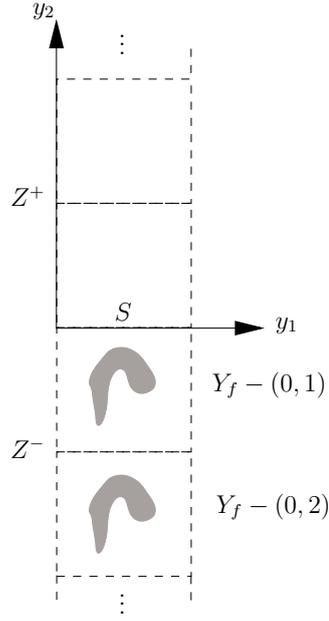}
\caption{The boundary layer geometry}\label{afotobl}
\end{figure}
\vskip1pt
 On the figure,
 the interface is $S=(0,1)\times \{ 0\} $, the free fluid slab is $Z^+ = (0,1) \times (0, +\infty )$
and the semi-infinite porous slab $Z^- =
\displaystyle\cup_{k=1}^\infty ( Y_F -\{ 0,k \} )$. The flow region
is then $Z_{BL} = Z^+ \cup S \cup Z^- $. \vskip3pt We consider the
following problem: \vskip1pt Find $\{ \beta^{j, bl} , \omega ^{j, bl} \} $, $j=1,2$,
with square-integrable gradients satisfying
\begin{subequations}\label{equ.boundary layer}
\begin{gather} -\Delta _y \beta^{j, bl} +\nabla_y \omega ^{j, bl} =0\qquad \hbox{ in
} Z^+ \cup Z^- \label{BJ4.2}\\ \; \mbox{ div } _y  \beta^{j, bl} =0\qquad \hbox{ in }
Z^+ \cup Z^- \label{4.3} \\ \bigl[ \beta^{j, bl} \bigr]_S (\cdot , 0)= K_{2j} \mathbf{e}^2 - \mathbf{w}^j
 \quad \hbox{ on } S \label{4.4}\\
 \bigr[ \{ \nabla_y \beta^{j, bl}
-\omega^{j, bl} I \}  \mathbf{e}^2 \bigl]_S (\cdot , 0) = - \{ \nabla_y \mathbf{w}^{ j}
-\pi^{ j} I \}  \mathbf{e}^2\ \hbox{ on }
S\label{4.5)} \\ \beta^{j, bl} =0 \quad \hbox{ on }
\displaystyle\cup_{k=1}^{\infty} ( \p Y_s -\{ 0,k \} ), \qquad \{
\beta^{j, bl} , \omega^{j, bl} \} \, \hbox{ is } 1- \hbox{periodic in }
y_1 .\label{4.6} \end{gather}
\end{subequations}
By Lax-Milgram's lemma, there exists  a
unique $\beta^{j, bl} \in L^2_{loc} (Z_{BL} )^2, \; \nabla _y \beta^{j, bl} \in L^2
(Z^+ \cup Z^- )^4$  satisfying (\ref{BJ4.2})-(\ref{4.6}) and
 $\omega^{j, bl} \in L^2_{loc}
(Z_{BL})$, which is unique up to a constant and satisfying
(\ref{BJ4.2}). After the results from \cite{JaMi2}, the  system (\ref{BJ4.2})-(\ref{4.6}) describes a
boundary layer, i.e.  $\beta^{j, bl} $ and $ \omega^{j, bl} $
stabilize exponentially towards  constants, when $\vert y_2\vert \to
\infty$: There exists $\gamma_0 >0$ and $\mathbf{C}^{j, bl}$ and $C^j_\pi$ such that
\begin{gather}
| \beta^{j, bl} - \mathbf{C}^{j, bl} | + | \omega^{j, bl} - C^j_\pi | \leq C e^{-\gamma_0 y_2}, \qquad y_2 >0,
\label{decay1} \\
    e^{-\gamma_0 y_2} \nabla_y \beta^{j, bl}  , \quad e^{-\gamma_0 y_2} \beta^{j, bl}, \quad  e^{-\gamma_0 y_2} \omega^{j, bl} \in L^2 (Z^-), \label{decay2} \\
    \mathbf{C}^{j, bl} = ({C}^{j, bl}_1 , 0) = (\int_S \beta^{j, bl}_1 (y_1 , +0) \ d y_1 , 0), \label{decay3} \\
C^j_\pi =\int^1_0 \omega^{j, bl} (y_1 , +0) \ d y_1 .\label{decay3A}
\end{gather}
The case $j=2$ is of special importance. If we suppose the mirror symmetry of the solid obstacle $Y_s$ with respect to $y_1$, then it is easy to prove that $w^2_1$ is uneven in $y_1$ with respect to the line $\{ y_1 = 1/2 \}$, and $w^2_2$ and $\pi^2$ are even. Consequently, $K_{12}=K_{21}=0$ and the permeability tensor $K$ is diagonal. Next we see that $\beta^{2, bl}_1$ is uneven in $y_1$ with respect to the line $\{ y_1 = 1/2 \}$, and $\beta^{2, bl}_2$ and $\omega^2$ are even. Using formula (\ref{decay3}) yields ${C}^{2, bl}_1 =0$ in the case of the mirror symmetry of the solid obstacle $Y_s$ with respect to $y_1$.

\subsection{The macroscopic model}

 Now we introduce the effective problem in $\Omega$.
 It consist of two problems, which are to be solved sequentially. The first problem is posed in $\Omega_2$ and reads:

Find a pressure field $P^D$ which
is the $L-$ periodic in $x_1$ function satisfying
\begin{subequations}\label{equ.darcy}
\begin{gather}
    \mathbf{u}^{D} = -K \nabla  P^D \quad \mbox{and} \quad  \mbox{ div } \bigg( K  \nabla  P^D \bigg) =0\; \mbox{ in } \; \O_2  \label{Presspm} \\
 P^D  =  0   \; \mbox{ on } \; \Sigma; \quad -K  \nabla  P^D  |_{\{ x_2 =-H \} } \cdot \mathbf{e}^2 =g.  \label{Presspm2}
\end{gather}
\end{subequations}

We note that the pressure field $P^D$ is equal to a particular constant, which is equal to zero.

Problem (\ref{Presspm})-(\ref{Presspm2}) admits a unique solution $\{ \mathbf{u}^{D} , P^D \} \in C^\infty ({\overline \Omega}_2)^3 $.
\vskip7pt
\noindent Next, we study the situation in the unconfined fluid domain  $\O_1$:

Find a velocity field
$\mathbf{u}^{eff}$ and a pressure field $p^{eff}$ such that
\begin{subequations}\label{equ.effective flow}
\begin{gather}
- \triangle \mathbf{u}^{eff}  + \nabla p^{eff} =0 \qquad \hbox{ in } \O_1 ,\label{4.91}\\
\; \mbox{ div }  \ \mathbf{u}^{eff} = 0 \qquad \hbox{ in } \O_1 ,  \qquad \int_{\O_1} p^{eff} \ dx =0,\label{4.92}\\ \mathbf{u}^{eff} =
\mathbf{v}^D  \qquad  \hbox{ on } (0,L) \times  \{h\}  
; \quad \mathbf{u}^{eff} \; \hbox{ and }  \;  p^{eff} \quad  \hbox{ are } \;
L-\hbox{periodic in} \quad x_1, \label{4.94}\\ u^{eff}_2 = -K  \nabla  P^D   \cdot \mathbf{e}^2 = - K_{22} \frac{\p P^D }{ \p
x_2 }  \qquad
\hbox{ and } \quad u^{eff}_1 = C^{2,bl}_1 \frac{\p P^D }{ \p
x_2 }  \quad \hbox{ on } \quad \Sigma  . \label{4.95}\end{gather}\end{subequations}
The constant $C^{2, bl}_1$ is  given by (\ref{decay1}) and requires solving problem (\ref{BJ4.2})-(\ref{4.6}).

Again, using the compatibility condition (\ref{Compa}) we obtain easily that problem (\ref{4.91})-(\ref{4.95}) has a unique solution $\{ \mathbf{u}^{eff} , p^{eff} \} \in C^\infty ({\overline \Omega}_1)^3 $.

\subsection{The main result}

In this section we formulate the approximated model. We expect that the Stokes system remains valid in $\Omega_1$.

Since $U_B =O(1) \neq 0$, the filtration velocity has to be of order $O(1)$.  Therefore, after \cite{All97}, \cite{JaegMik09} and \cite{SP80}, the asymptotic behavior of the velocity and pressure fields in the porous part $\Omega_2$, in the limit $\ep \to 0$, is given by the two-scale expansion
\begin{gather*}
 \mathbf{v}^\ep  \approx   \mathbf{u}^{0} (x , y ) + \ep  \mathbf{u}^{1} (x , y ) + O(\ep^2 ), \quad y=\frac{x}{\ep },\\
  p^\ep \approx \frac{1}{\ep^2} P^D (x) + \frac{1}{\ep} p^1 (x,y) + O(1 ), \quad y=\frac{x}{\ep }, \\
\mathbf{u}^{0} (x , y ) = - \sum_{j=1}^2 \mathbf{w}^j (y) \frac{\p P^D (x)}{\p x_j}, \quad p^1 (x,y) = - \sum_{j=1}^2 \pi^j (y) \frac{\p P^D (x)}{\p x_j}.
\end{gather*}
The boundary layers given by (\ref{BJ4.2})-(\ref{4.6}) will be used to link the above approximation on $\Sigma$ with the solution of the Stokes system. With such strategy,
at the main order approximation  reads
\begin{gather}
\mathbf{v}^\ep = H(x_2) (\mathbf{u}^{eff}  -  C^{2, bl}_1   \frac{\p P^D }{ \p
x_2 }  |_\Sigma \mathbf{e}^1 ) - H(-x_2) \sum_{k=1}^2 \frac{\p P^D }{ \p
x_k} \mathbf{w}^k (\frac{x}{\ep})  +\notag \\
   \frac{\p P^D }{ \p
x_2} |_\Sigma \beta^{2, bl} (\frac{x}{\ep}) + O(\ep)  +\; \mbox{outer  boundary layer},
\label{Veloc1} \\
p^\ep =
 H(x_2) p^{eff} + H(-x_2 ) \{ \ep^{-2} P^D - \frac{1}{\ep}  \sum_{k=1}^2 (\frac{\p P^D } {\p
x_k}  \pi^k (\frac{x}{\ep}) + \mathcal{A}^k_\pi \delta_{2k}) \}  +\notag \\
  \frac{1}{\ep }  (\omega^{2, bl} (\frac{x}{\ep})   - C^2_\pi H(x_2 ) ) \frac{\p P^D }{ \p
x_2} |_\Sigma +  o(\frac{1}{\ep}) +\; \mbox{outer  boundary layers},
 \label{press1}
\end{gather}
where  $H(t)$ is the Heaviside function.
We will see that $\displaystyle \mathcal{A}^k_\pi  = C^2_\pi \frac{\p P^D }{ \p
x_2} |_\Sigma $.


\begin{theorem} \label{Thbasic}
Let $\mathcal{O}$ be a neighborhood of $x_2 =-H$.
Let us suppose the geometry and data smoothness as above and the compatibility condition (\ref{Compa}). Let  $p^\ep$ be extended to $\Omega_2$ by formula (\ref{Extpress}) of Lipton and Avellaneda.
 Then   we have
\begin{gather}
|| \mathbf{v}^\ep - \mathbf{u}^{eff} ||_{L^2 (\Omega_1)} \leq C\sqrt{\ep}
\label{Est1} \\
|  || \mathbf{v}^\ep +    \frac{\p P^D }{ \p
x_2 }  |_\Sigma  ( K_{22} \mathbf{e}^2 -  \beta^{2, bl} (\frac{x_1}{\ep} , 0+)) ||_{L^2 (\Sigma )} \leq C\sqrt{\ep}
\label{Est2} \\
|| \mathbf{v}^\ep + \sum_{k=1}^2 \frac{\p P^D }{ \p
x_k} \mathbf{w}^k (\frac{x}{\ep}) -   \frac{\p P^D }{ \p
x_2} |_\Sigma \beta^{2, bl} (\frac{x}{\ep}) ||_{L^2 (\Omega_2  \setminus \mathcal{O})} \leq C \ep
\label{Est3}\\
  || p^\ep - H(-x_2 ) \ep^{-2} P^D ||_{L^2 (\Omega )} \leq \frac{C}{\ep} . \label{Est4}
\end{gather}
\end{theorem}

Inspection of the proof of theorem (\ref{Thbasic})   shows that we can obtain slightly better estimates by rearranging the term $$ \frac{1}{\ep }  (\omega^{2, bl} (\frac{x}{\ep})   - C^2_\pi H(x_2 ) ) \frac{\p P^D }{ \p
x_2} |_\Sigma .$$

We obtain
\begin{theorem} \label{Thadvanced}
Let $\mathcal{O}$ be a neighborhood of $x_2 =-H$.
Let us suppose the geometry and data smoothness as above and the compatibility condition (\ref{Compa}). Let  $p^\ep$ be extended to $\Omega_2$ by formula (\ref{Extpress}) of Lipton and Avellaneda.
 Then   we have
\begin{gather}
|| \mathbf{v}^\ep - \mathbf{u}^{eff} + C^{2, bl}_1   \frac{\p P^D }{ \p
x_2 }  |_\Sigma \mathbf{e}^1 - \frac{\p P^D }{ \p
x_2} |_\Sigma \beta^{2, bl} (\frac{x}{\ep}) ||_{L^2 (\Omega_1)} \leq C\ep
\label{Est1A} \\
|| \mathbf{v}^\ep +    \frac{\p P^D }{ \p
x_2 }  |_\Sigma  ( \mathbf{w}^k (\frac{x_1 }{\ep} , 0-)  -  \beta^{2, bl} (\frac{x_1}{\ep} , 0-)) ||_{L^2 (\Sigma )} =\notag \\
  || \mathbf{v}^\ep +    \frac{\p P^D }{ \p
x_2 }  |_\Sigma  ( K_{22} \mathbf{e}^2 -  \beta^{2, bl} (\frac{x_1}{\ep} , 0+)) ||_{L^2 (\Sigma )} \leq C\ep
\label{Est2A} \\
|| \mathbf{v}^\ep + \sum_{k=1}^2 \frac{\p P^D }{ \p
x_k} \mathbf{w}^k (\frac{x}{\ep}) -   \frac{\p P^D }{ \p
x_2} |_\Sigma \beta^{2, bl} (\frac{x}{\ep}) ||_{L^2 (\Omega_2  \setminus \mathcal{O})} \leq C \ep
\label{Est3A}\\
  || p^\ep - H(-x_2 ) (\ep^{-2} P^D - \ep^{-1} ( C^2_\pi  \frac{\p P^D }{ \p
x_2} |_\Sigma ) +\sum_{j=1}^2 \pi^{j} (\frac{x}{\ep}) \frac{\p P^D}{\p x_j}))||_{L^2 (\Omega )} \leq \frac{C}{\sqrt{\ep}} . \label{Est4A}
\end{gather}
\end{theorem}

\begin{remark} We took as correction to $P^D$ the quantity $-C^2_\pi \partial_{x_2} P^D |_\Sigma$. In fact the better choice would be to take a function satisfying equations (\ref{Presspm})-(\ref{Presspm2}), with value on $\Sigma$ being  $-C^2_\pi \partial_{x_2} P^D |_\Sigma$, instead of zero. Since the order of approximation does not change, we make the simplest possible choice.

If the effective porous medium pressure is $P^{D, eff} = P^D -\ep C^2_\pi \partial_{x_2} P^D |_\Sigma$, then the requirement that we can only have an $O(1)$ normal stress jump on $\Sigma$ yields
\begin{equation}\label{PressureSlip}
  P^{D, eff} + \ep C^2_\pi \partial_{x_2} P^{D, eff} = O(\ep^2) \quad \mbox{on} \quad \Sigma.
\end{equation}
The relation (\ref{PressureSlip}) indicates presence of an effective pressure slip at the interface $\Sigma$.
 Since $\ep$ is related to the square root of the permeability, in the dimensional formulation, our result compares to the numerical experiments by Sahraoui and Kaviany in \cite{K95} and \cite{SahKav92}. They found it being small for parallel flows. We find it small but of order of the corrections in the law by Beavers and Joseph in the case of the transverse flow.
\end{remark}

\begin{remark} We obtain an expression for the tangential slip velocity. Since it is zero in the isotropic case, we do not confirm formulas like (86), page 2645, from \cite{OTW1:95} or like formula  (31) for oblique flows from \cite{SahKav92}, which generalize the law by Beavers and Joseph.

In  \cite{OTW1:95}, formula (71), page 2643, expresses the continuity of the averaged velocities. By construction, we have the trace continuity for our approximation. Nevertheless, one usually does not keep the boundary layers in the macroscopic model.
If we eliminate the boundary layers and all low order terms, the tangential effective velocity at the interface $\Sigma$ is
$$ \mathbf{u}^{eff} =C^{2,bl}_1 \frac{\p P^D }{ \p
x_2 }  \mathbf{e}^1  - K_{22} \frac{\p P^D }{ \p
x_2 } \mathbf{e}^2  , $$
(see (\ref{4.95}) and from the porous media side
$$ \mathbf{u}^{D} = - K_{22} \frac{\p P^D }{ \p
x_2 } \mathbf{e}^2  . $$
Therefore we find out that there is an effective tangential velocity jump at the interface.
\end{remark}

\section{Numerical confirmation of the effective interface conditions}\label{numerics}
\label{Numerics}
This section is dedicated to the numerical confirmation of the analytical results shown above.
We solve the problems needed to numerically compare the microscopic with the macroscopic problem by the finite element method (FEM). For the FEM theory we refere to standard literature, e.g., \cite{Ciarlet:2002} or \cite{BrennS:2002}.

For the discretization of the Stokes system we use the Taylor-Hood element, which is inf-sup stable \cite{BrezzF:1991}, therefore it does not need stabilization terms.
In particular, since the homogenization error in some of the proposed estimates is small in comparison with the discretization error even for meshes with a number of elements in the order of millions, we have used higher order finite elements (polynomial of third degree for the velocity components and of second degree for the pressure) to reduce the discretization error.

\begin{figure}
 \centering
\begin{tabular}{cc}
 \subfloat[Circle]{\resizebox{0.3\textwidth}{!}{\includegraphics[trim=54mm 80mm 21mm 83mm, clip, width =0.8\textwidth]{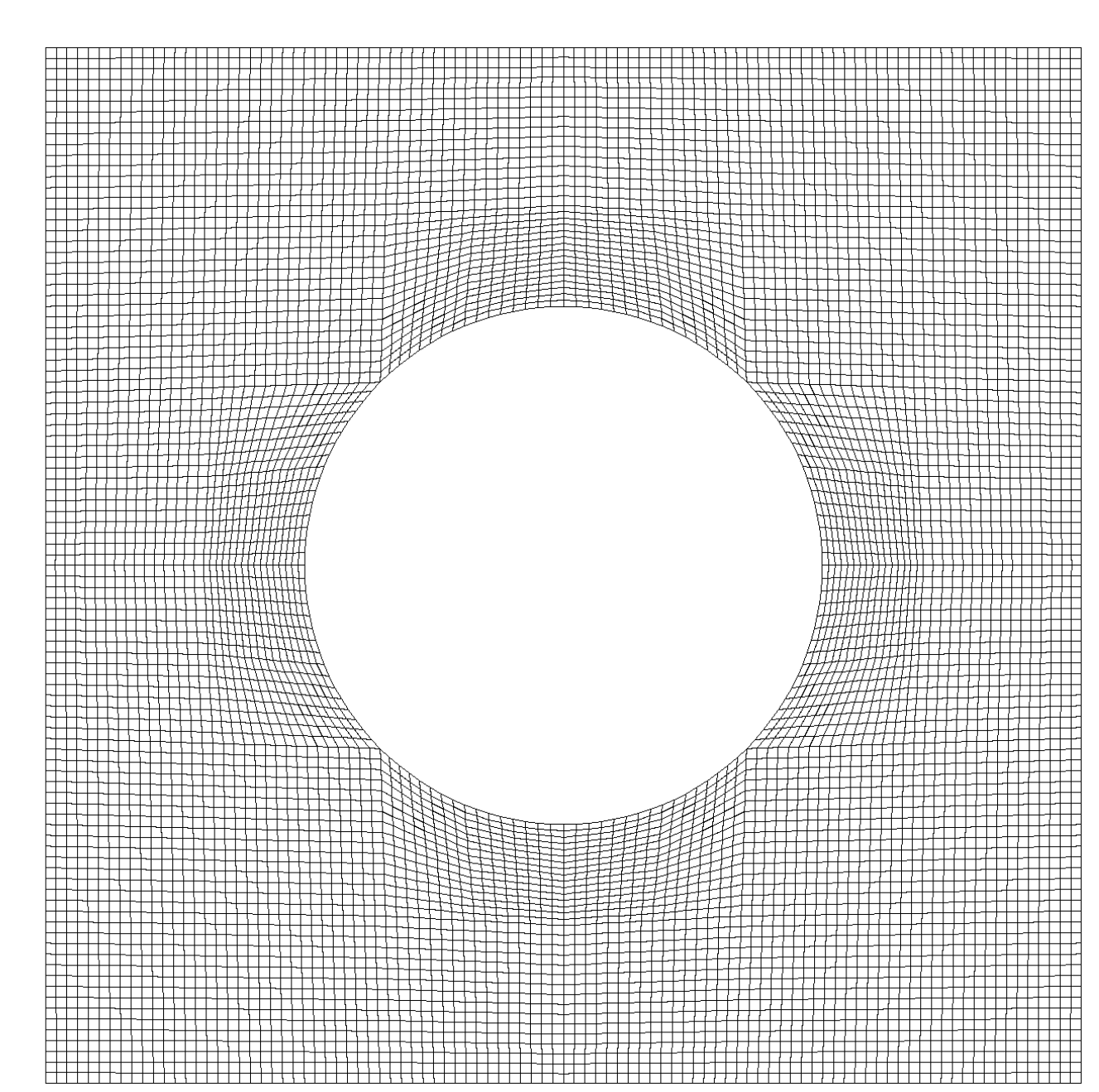}}}
 \hspace{0.1\textwidth}
 \subfloat[Ellipse]{\resizebox{0.3\textwidth}{!}{\includegraphics[trim=54mm 80mm 21mm 83mm, clip,width=0.8\textwidth]{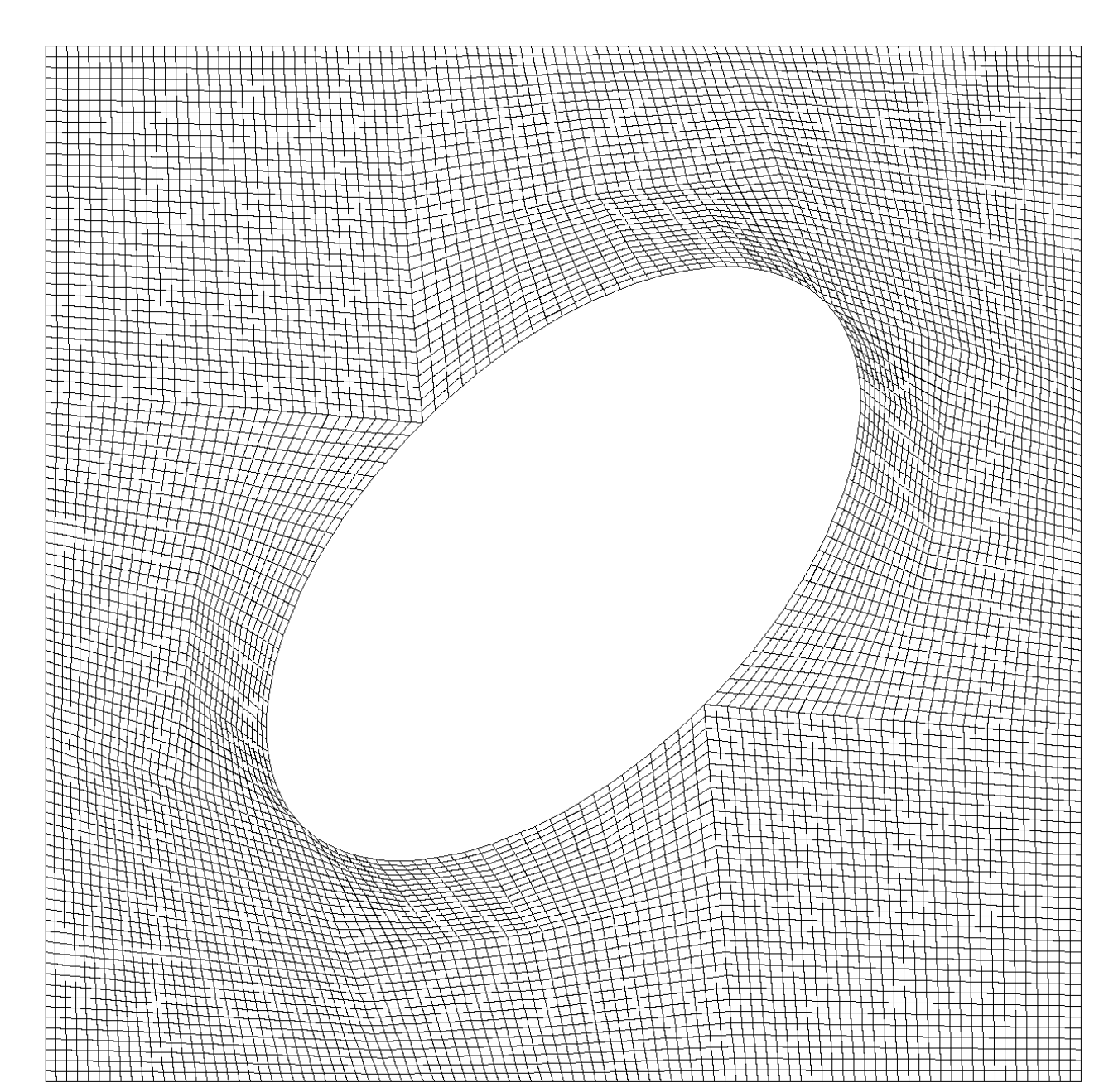}}}
\end{tabular}
\caption{Mesh of the fluid part of the unit cell for the two types of inclusions: circles and ellipses}\label{fig.inclusions}
 \end{figure}

The flow properties depend on the geometry of the pores. In particular there is a substantial difference between the case with symmetric inclusions with respect to the axis orthogonal to the interface and the case with asymmetric inclusions.
We use therefore two different types of inclusion in the porous part, {\it circles} and rotated {\it ellipses}, i.e. ellipses with the major principal axis non parallel to the flow.
The increased accuracy using higher order finite elements in the numerical solutions was necessary, as shown later, especially for the case with symmetric inclusions.
The geometries of the unit cells $Y = (0.1)^2$, see figure~\ref{fig.inclusions}, for these two cases are as follows:
\begin{enumerate}
\item the solid part of the cell $Y_s$ is formed by a circle with radius $0.25$ and center $(0.5, 0.5)$.
\item $Y_s$ consists of an ellipse with center $(0.5, 0.5)$ and semi-axes $a=0.4$ and $b=0.2$, which are rotated anti-clockwise by $45^\circ$.
\end{enumerate}
In addition, since the considered domains have curved boundaries we use cells of the FEM mesh with curved boundaries (a mapping with polynomial of second degree was used for the geometry) to obtain a better approximation.

All computations are done using the toolkit \texttt{DOpElib} (\cite{GollWW:2012}) based upon the C++-library \texttt{deal.II} (\cite{BangeHK:2007}).
\subsection{Numerical setting}\label{sec:numerical setting}
In this subsection we describe the setting for the numerical test. To confirm the estimates of Theorem \ref{Thbasic} and \ref{Thadvanced} we have to solve the microscopic problem \eqref{equ.micro} to get $\ve$ and \pe{}, the macroscopic problems \eqref{equ.darcy} and \eqref{equ.effective flow} to get $\mathbf{u}^{eff}, p^{eff}$ and $P^D$, the cell problem \eqref{1.57Cell} to calculate the permeability tensor $K$, the velocity vector $\mathbf{w}$ and pressure $\pi$, and the boundary layer \eqref{equ.boundary layer} for the velocity $\boldsymbol{\beta}^{bl}$ and pressure $\omega^{bl}$.

To reduce the discretization errors we consider a test case, described below, for which it is easy to derive the exact form of the macroscopic solution. As we will show below, the analytical solution of the macroscopic problem can be expressed in terms of the solution of the cell and boundary layer problems. The discretization error of the macroscopic problem in this case depends on the discretization error of the cell and boundary layer problems and does not imply therefore an additional discretization error.

We consider the following domains $\Omega = [0,1]\times[-1,1]$ and $\Oe = \Omega \setminus \text{`the obstacles'}$, where the obstacles are either circles or ellipses as described in the subsection above. In our example we consider the in- and outflow condition
\begin{align}
\mathbf v_D &= (0,-1) \text{ and } g =-1.
\end{align}
in the microscopic problem \eqref{equ.micro}

The macroscopic solution in this setting is
\begin{subequations}
\label{macroscopic solution}
\begin{align}
\label{u1eff} u^{\eff}_1 &= \frac {\cbl}{K_{22}} (1-y),\\
\label{u2eff}  u^{\eff}_2 &= -1,\\
\label{peff} p^{\eff}&=0,\\
\label{P0} \pnull &= \frac 1 {K_{22}} y.
\end{align}
\end{subequations}

The macroscopic solution depends on the solution of the cell problem though the permeability $K$, see expressions \eqref{u1eff} and \eqref{P0}. Furthermore it depends on the solution of the boundary layer though the constant $\cbl$. The macroscopic problems \eqref{equ.darcy} and \eqref{equ.effective flow} are therefore not numerically solved.

The microscopic problem \eqref{equ.micro} is solved with around 10--15 million degrees of freedom, the cell problem uses around 7 million degrees of freedom. The permeability constant has been precisely calculated using the goal oriented strategy for mesh adaptivity described in \cite{CGMCM:13}.

In the boundary layer problem, due to the interface condition \eqref{4.4}, the velocity as well as the pressure are discontinuous on the interface $S$. Since with the $H^1$ conform finite elements chosen for the discretization the discontinuity cannot be properly approximated, we have decided to transform the problem so that the solution variables are continuous across $S$. The values of $\boldsymbol \beta^{bl}$ and $\pi^{bl}$ needed to check the estimates are recovered by post-processing.
For the numerical solution, as explained in detail in the appendix of our previous work \cite{CGMCM:13}, we use a cut-off domain, which is justified by the exponential decaying of the boundary layer solution.
The solution of the boundary layer problem is obtained with a mesh of around 4 million degrees of freedom and the constants \cbl{} and \cblpi{} are calculated by the goal oriented strategy for mesh adaptivity described in \cite{CGMCM:13} where we have made sure that the cut-off error is smaller than the discretization error. We note that in the computation of  $C^2_\pi$ we do not use the formula given in  \eqref{decay3A} but the equivalent one
\begin{align}
\int^1_0 \omega^{j, bl} (y_1 , 1) \ d y_1
\end{align}
as this proved to be advantageous numerically.

In table~\ref{tab.results constants} the computed constants $K, \cbl$ and $\cblpi$ for the two different inclusions are listed. As the permeability tensor $K$ has for the given cases the form $$K=\begin{pmatrix}K_{11} &K_{12}\\K_{12} &K_{11}\end{pmatrix},$$ i.e. it holds $K_{11}=K_{22}$ and $K_{12}=K_{21}$, we state only $K_{11}$ and $K_{12}$. Additionally, we give an estimation of the discretization error.
\begin{table}
\centering
\begin{tabular}{c|ll}
\toprule
&\textbf{circular inclusions}& \textbf{oval inclusions}\\
\midrule
$K_{11}$  	& 0.0199014353519271 $\pm2\cdot10^{-12}$	&0.0122773324576884 $\pm2\cdot10^{-13}$	\\
$K_{12}$  	& 0 					&0.00268891986291451 $\pm2\cdot10^{-13}$	\\
\cbl  		&0						&-0.003336740001686 $\pm4\cdot10^{-10}$	\\
\cblpi 		&0.025777570627281 $\pm3\cdot10^{-8}$		&-0.004429782196436	$\pm1\cdot10^{-8}$\\	
\bottomrule
\end{tabular}
\caption{Computed constants for the computations in the example.}\label{tab.results constants}
\end{table}

\subsection{Numerical results}\label{sec:numerical results}
In this section we present the numerical confirmation of the convergence rates of the homogenization errors (\ref{Est1}--\ref{Est4}) and (\ref{Est1A}--\ref{Est4A}).
 \begin{figure}[!htb]
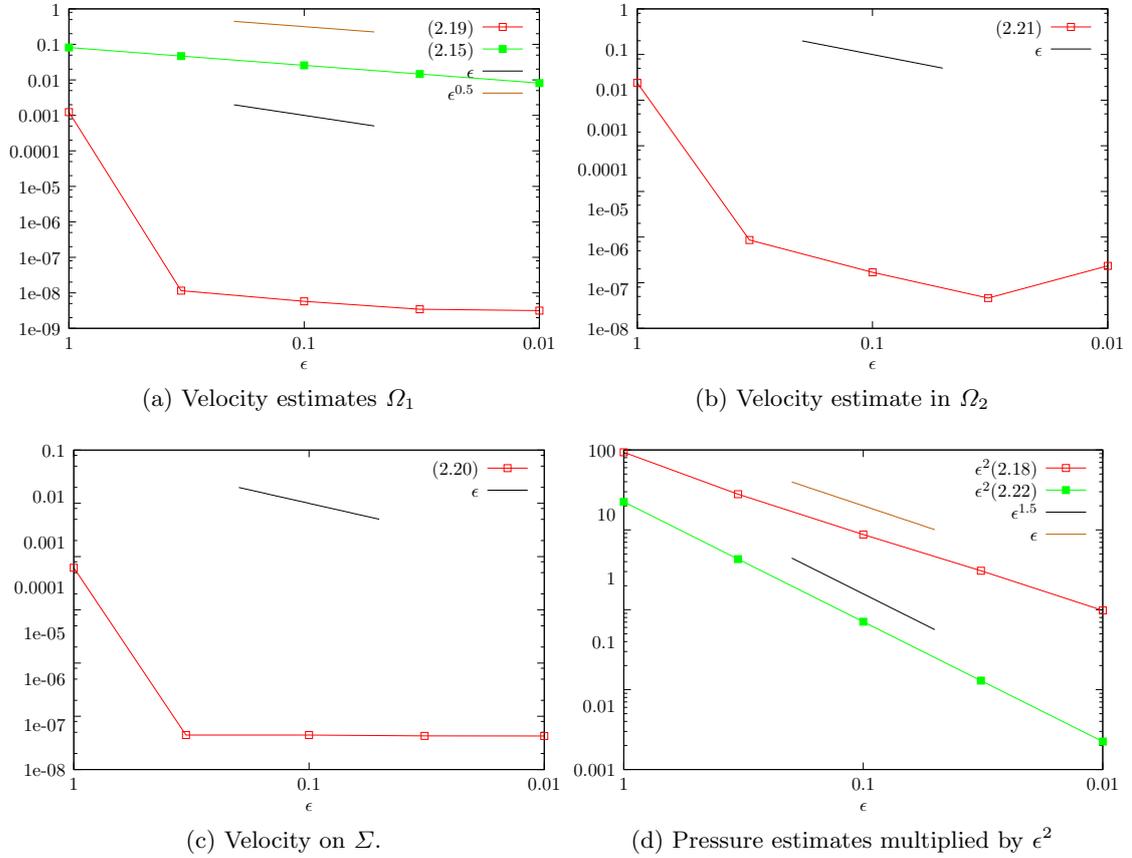

\centering
\begin{tabular}{cc}
 \subfloat[Velocity estimates $\Omega_1$]{\resizebox{0.5\textwidth}{!}{\input{figures/ConvergenceCircleFluid32.tex}}\label{fig.convergence circs fluid}}
 \subfloat[Velocity estimate in $\Omega_2$]{
 \label{fig.convergence circs omega_2}
 \resizebox{0.5\textwidth}{!}{\input{figures/ConvergenceCirclePorous32.tex}}}
 \\
 \subfloat[Velocity on $\Sigma$.]{\label{fig.convergence circs sigma}\resizebox{0.5\textwidth}{!}{\input{figures/ConvergenceCircleFace32.tex}}}
 \subfloat[Pressure estimates multiplied by $\epsilon^2$]{\resizebox{0.5\textwidth}{!}{\input{figures/ConvergenceCirclePressure32.tex}}\label{fig.convergence circs pressure}}
\end{tabular}
\caption{Convergence results for circles.}\label{fig.convergence circs}
 \end{figure}

 For our test we set $\Omega_2\setminus \mathcal O = [0,1]\times[-0.6,0]$, and we use a computation of the boundary layer on a cut-off domain ranging from $-4$ to $4$. This means that to compute the norms we evaluate the terms involving the boundary layer only for $x\in \Omega$ with $-4\epsilon < x_2<4 \epsilon$. Outside of this region we assume the difference between the boundary layer components and their respective asymptotic values to be sufficiently small.

In the case of inclusions symmetric in the sense explained above, e.g. circles, the homogenization errors are much smaller than the numerical error even for large epsilon such as $0.1$ as can be observed in figure~\ref{fig.convergence circs}).
The lines with markers represent the results of the computations for  $\epsilon \in \{1, \tfrac 13, 0.1, \tfrac 1 {31}, 0.01\}$, the solid lines are reference values for various convergence rates and are plotted only to compare the respective slopes.

\begin{figure}[!htb]
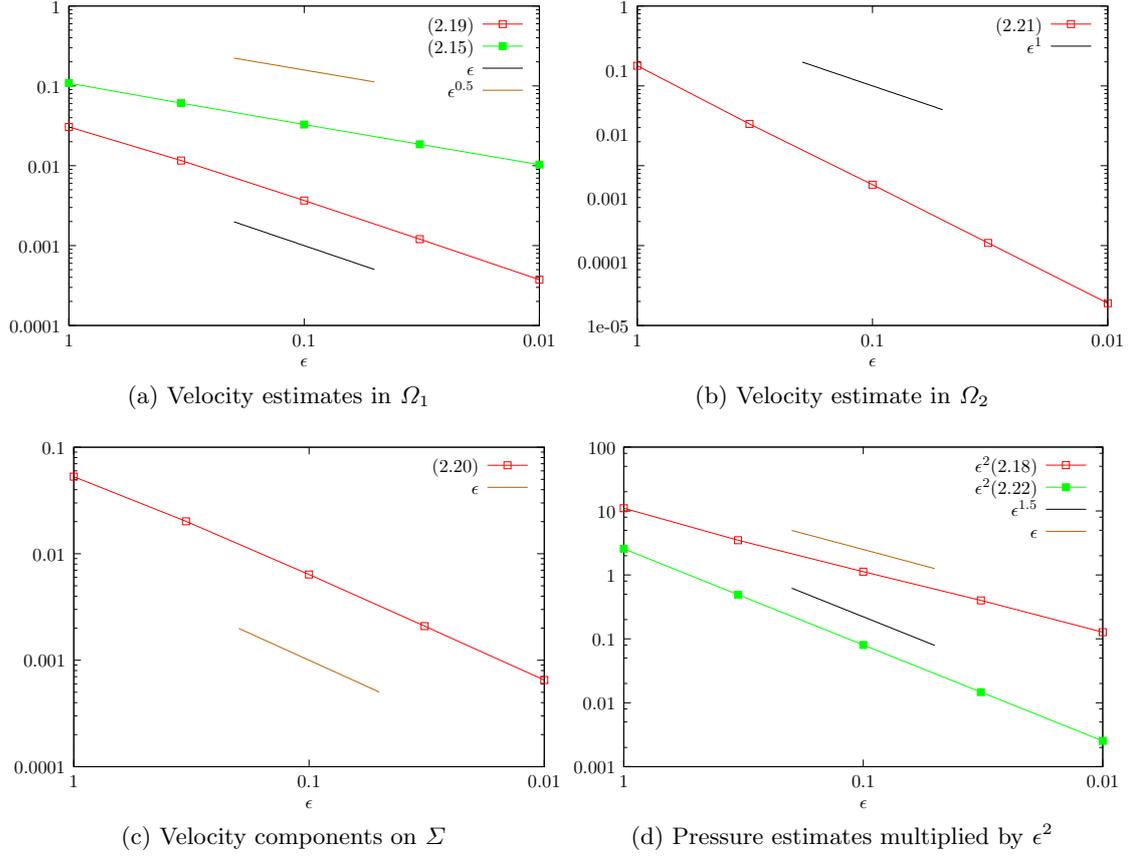

\centering
\begin{tabular}{cc}
 \subfloat[Velocity estimates in $\Omega_1$]{\resizebox{0.5\textwidth}{!}{\input{figures/ConvergenceEllipseFluid32.tex}}\label{fig.convergence ellipse fluid}}
 \subfloat[Velocity estimate in $\Omega_2$]{
 \resizebox{0.5\textwidth}{!}{\input{figures/ConvergenceEllipsePorous32.tex}}\label{fig.convergence ellipse porous}}
 \\
 \subfloat[Velocity components on $\Sigma$]{\resizebox{0.5\textwidth}{!}{\input{figures/ConvergenceEllipseFace32.tex}}\label{fig.convergence ellipse face}}
 \subfloat[Pressure estimates multiplied by $\epsilon^2$]{\resizebox{0.5\textwidth}{!}{\input{figures/ConvergenceEllipsePressure32.tex}}\label{fig.convergence ellipse pressure}}
\end{tabular}
\caption{Convergence results for elliptical inclusions.}\label{fig.convergence ell}
 \end{figure}

The case of circles is shown in figure \ref{fig.convergence circs}.
For the velocity in the fluid part of the domain the estimate \eqref{Est1} can be verified.
For the better estimate \eqref{Est1A}, that uses correction terms to improve the estimation, the homogenization error is so small that the curve shows only the numerical error, that in our case is only due to the discretization error since the quadrature error and the tolerance of the solver are smaller. In Figure \ref{fig.convergence circs omega_2} we can confirm \eqref{Est3A} only for values of epsilon not bigger than $\frac{1}{31}$, for $\epsilon =0.01$ the numerical error dominates the homogenization error.
In the estimates for circles  on the interface (Figure \ref{fig.convergence circs sigma}) we can observe only the numerical error for the same reason explained above.
Notice that the error for circles shown in Figure \ref{fig.convergence circs} is much smaller than the error for ellipses shown in Figure \ref{fig.convergence ell}.
In addition, we could verify both estimates for the pressure \eqref{Est4} and \eqref{Est4A} as shown in Figure \ref{fig.convergence circs pressure}. Note that the pressure estimates have been scaled multiplying by $\epsilon^2$.

The case of ellipses is shown in figure \ref{fig.convergence ell}.
As it can be observed, all estimates could be numerically verified, since the discretization error in this case was smaller than the homogenization error. Also in this case the pressure estimates have been scaled multiplied by $\epsilon^2$. Note, that we observe for the velocity in the porous domain a convergence rate of 1.5 instead of the predicted first order convergence, see figure~\ref{fig.convergence ellipse porous}.

\begin{figure}[htb]
 \subfloat[Streamlines of $\ve$\label{subfig.vc}]{\resizebox{0.21\textwidth}{!}{\includegraphics[trim=-20mm 150mm 10mm 15mm,clip,]{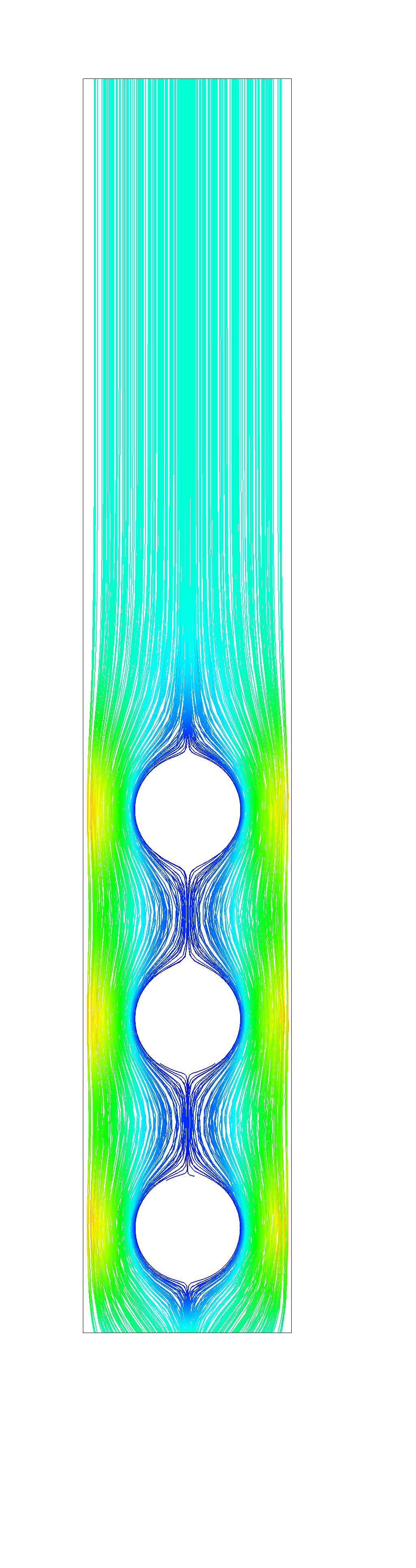}}}
 \subfloat[$\pe$\label{subfig.pc}]{\resizebox{0.3\textwidth}{!}{\includegraphics[trim=0mm 60mm 130mm 15mm]{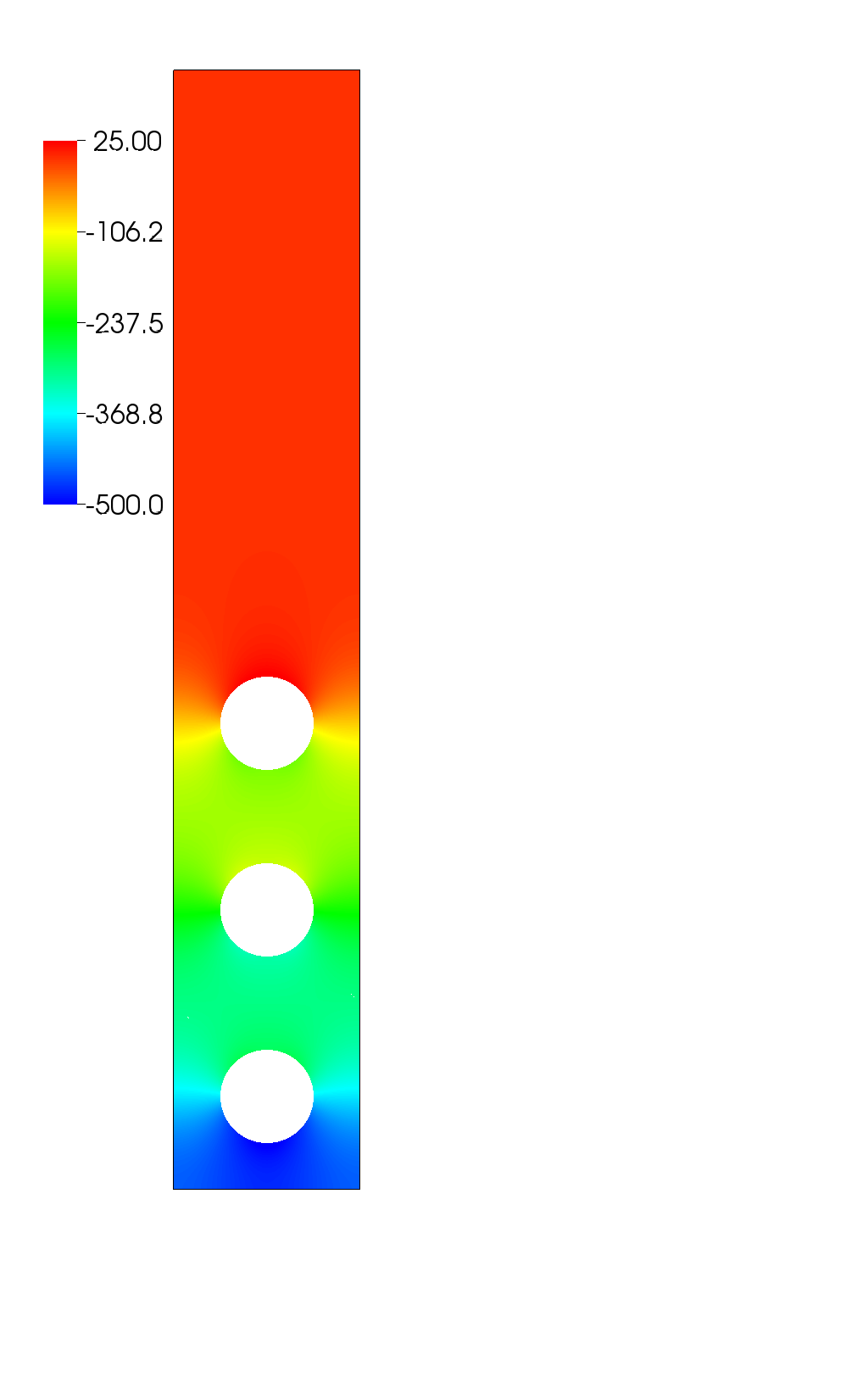}}}
 \subfloat[Streamlines of $\ve$\label{subfig.ve}]{\resizebox{0.21\textwidth}{!}{\includegraphics[trim=-20mm 150mm 10mm 15mm]{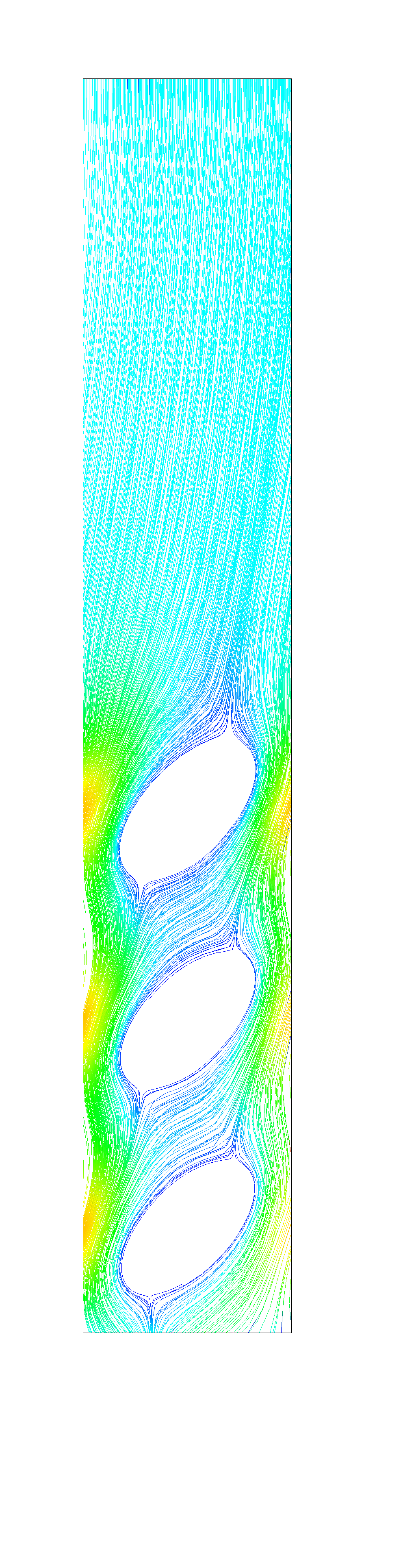}}}
 \subfloat[$\pe$\label{subfig.pe}]{\resizebox{0.3\textwidth}{!}{\includegraphics[trim=0mm 60mm 130mm 15mm]{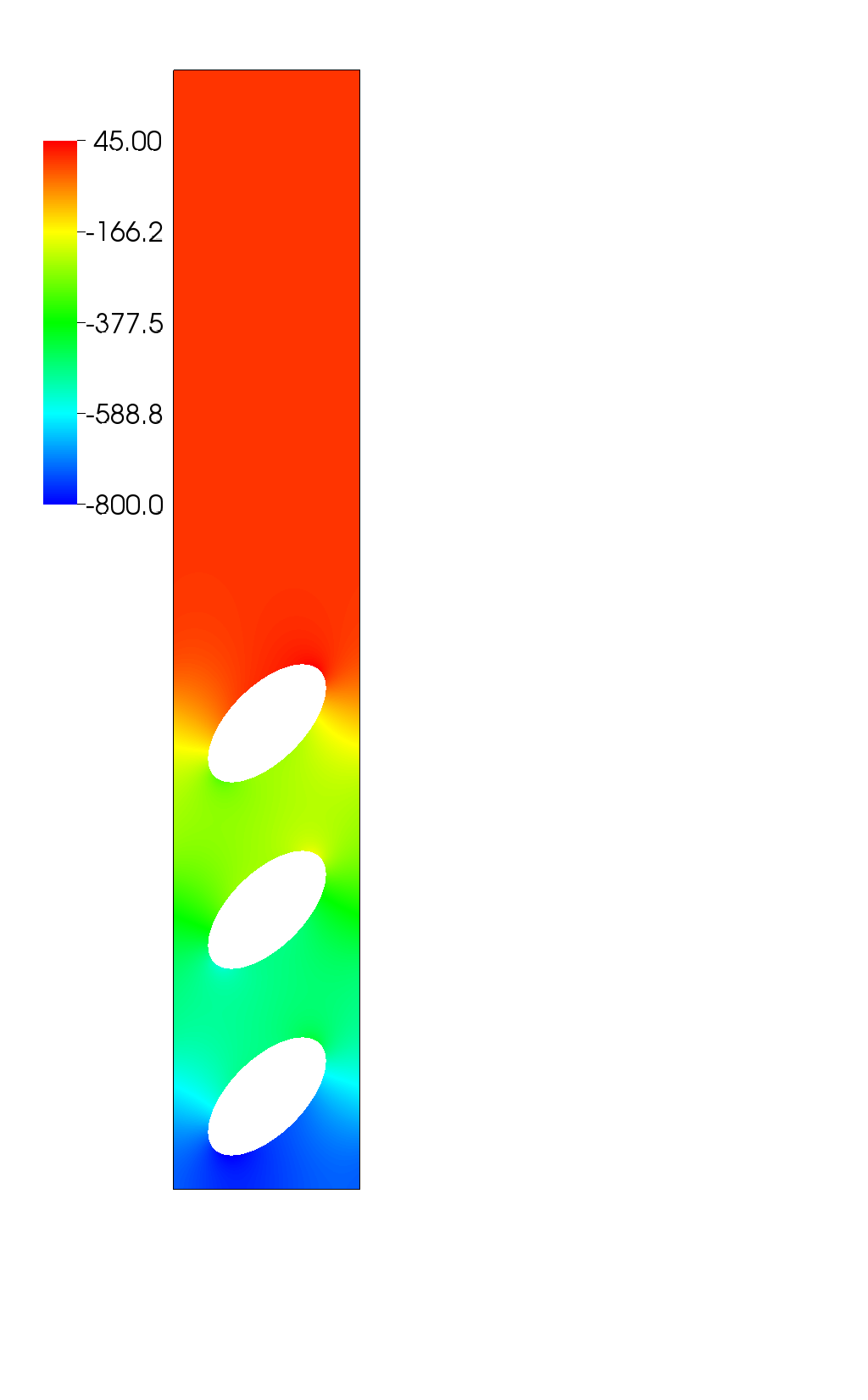}}}
\caption{Visualization of the solution to the microscopic problem for $\epsilon = \tfrac 13$. Subfigures \protect\subref{subfig.vc}and \protect\subref{subfig.pc} show the results for circular inclusions,  \protect\subref{subfig.ve}and \protect\subref{subfig.pe} for elliptical inclusions.}\label{fig.micro solution}
 \end{figure}

\begin{figure}[htb]
%
%
%
  \subfloat[Circular inclusions.\label{subfig.micro sol 3dv1 circ}]{\resizebox{0.5\textwidth}{!}{\includegraphics[trim=10mm 10mm 5mm 15mm, clip]{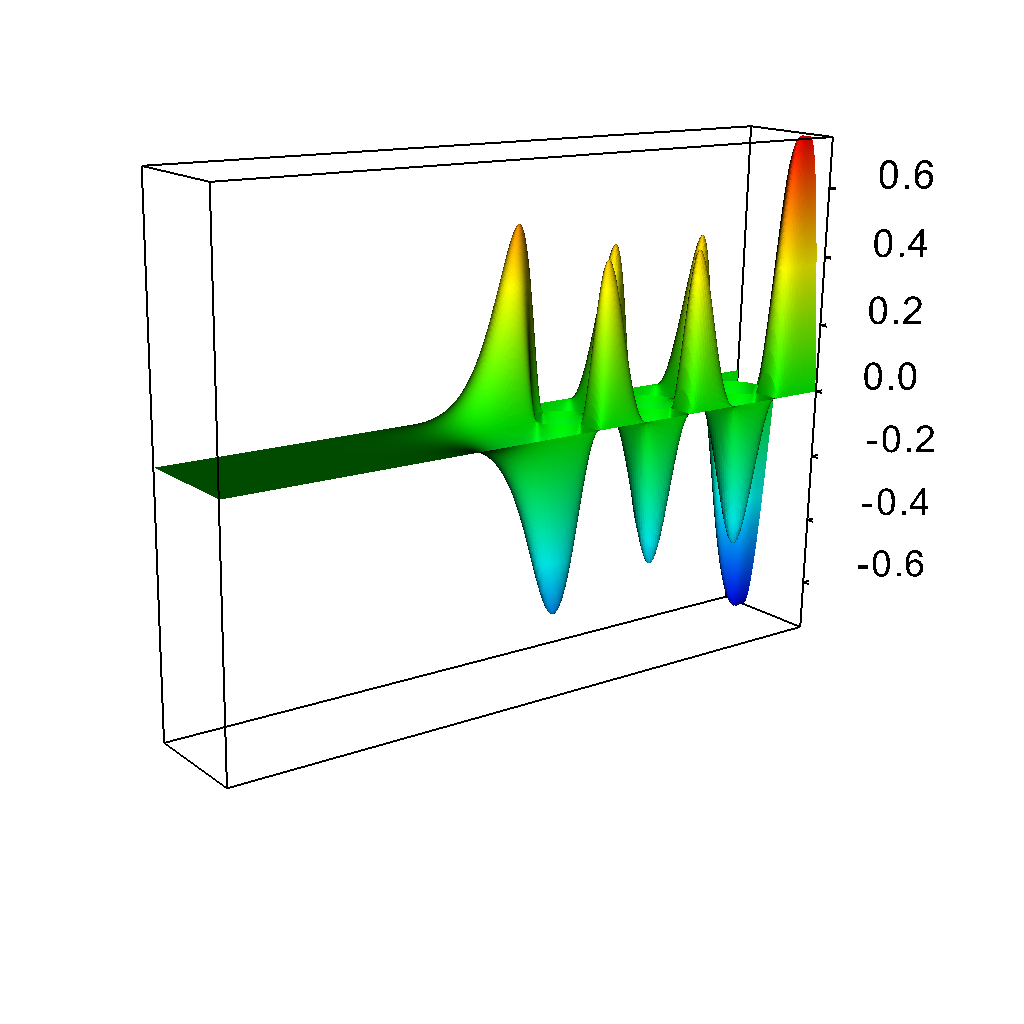}}}
 \subfloat[Elliptical inclusions.\label{subfig.micro sol 3dv1 ell}]{\resizebox{0.5\textwidth}{!}{\includegraphics[trim=10mm 10mm 5mm 15mm, clip]{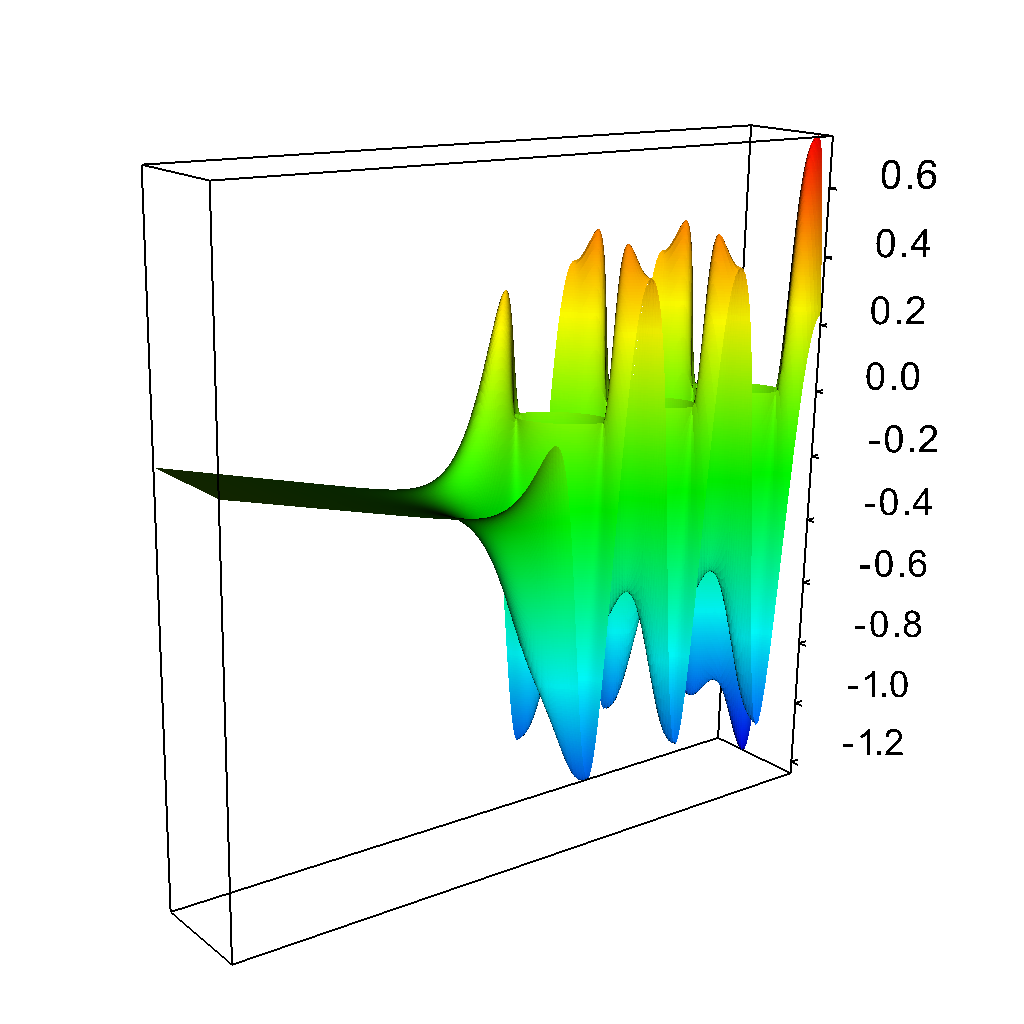}}}
\caption{Visualization of $\vecomp{1}$ for $\epsilon = \tfrac 13$.}\label{fig.micro solution 3dv1}
 \end{figure}

In conclusion, we show in figure \ref{fig.micro solution} and figure \ref{fig.micro solution 3dv1} pictures of the flow for the case $\epsilon=\tfrac 1 3$. Since we use periodic boundary conditions in the $x_1$-direction, constant in- and outflow data as well as a periodic geometry, the computations have been performed on a stripe of one column of inclusions to reduce the computational effort. In figure~\ref{subfig.vc} and \ref{subfig.ve} we see streamline plots of the velocity, figure~\ref{subfig.pc} and \ref{subfig.pe} show the corresponding pressures. Both pressures are nearly constant in the fluid part and show then a linear descent to the outflow boundary, similar to the effective pressure \eqref{peff} and \eqref{P0}.

Figure \ref{fig.micro solution 3dv1} shows only the values of the tangential velocity component. In the case of circular inclusions (figure~\ref{subfig.micro sol 3dv1 circ}), the velocity is nearly  zero throughout the fluid region and shows some oscillations around the mean value zero on the position of the interface. Note that the effective model prescribes here a no slip condition because it holds  $\cbl=0$. In figure~\ref{subfig.micro sol 3dv1 ell}) we see the corresponding solution for oval inclusions. We notice a linear descent from the inflow boundary (which lies in this picture on the left hand side) to the interface, which leads to the slip condition for the tangential velocity component of the effective flow in this case. Both behaviors are predicted from the effective interface condition for this velocity component, see \eqref{4.95}.


\section{Proof of Theorem \ref{Thbasic} via incremental accuracy correction}\label{Proofs1}

In the proofs which follow we will frequently use the space
\begin{gather}
\hskip-15pt V_{per} (\Oe ) = \{ \mathbf{z}\in H^1 (\Oe )^2 :
\, \mathbf{z} =0 \, \hbox{ on } \p \Oe \setminus \p \O , \;
\mathbf{z} =0 \quad \mbox{on} \quad \{ x_2 =h\}, \notag \\
 z_2 =0 \quad \mbox{on} \quad \{ x_2 =-H\} \quad
\hbox{ and } \mathbf{z} \, \hbox{ is L-periodic in } x_1 \hbox{ variable } \} .
\label{2.1} \end{gather}
We will follow the strategy from \cite{JaMi2}, write a variational equation for the errors in velocity and in pressure and reduce the forcing term in several steps.
We will frequently use the notation
\begin{equation}\label{1.57}
  \mathbf{w}^{j, \ep} (x) = \mathbf{w}^{j} (\frac{x}{ \ep}) \quad \mbox{and} \quad \pi^{j, \ep}(x) = \pi^{j} (\frac{x}{ \ep}),
\end{equation}
where $\{ \mathbf{w}^{j} , \pi^j \}$ is given by (\ref{1.57Cell}).
\vskip12pt
\subsection{Incremental accuracy correction, the 1st part}
\begin{proposition}\label{Prop1}
Let $P^D$  given by (\ref{Presspm})-(\ref{Presspm2}), $\{ \mathbf{u}^{eff} , p^{eff} \}$ be the solution for (\ref{4.91})-(\ref{4.95})
and $\{ \mathbf{w}^{j, \ep} , \pi^{j, \ep} \} $  defined by (\ref{1.57}). Let $\{ \mathbf{v}^\ep, p^{\ep} \}$ be the solution for (\ref{1.3})-(\ref{1.8}). Then
for every $\varphi \in V_{per} (\Oe )$ we have
\begin{gather}
| \langle \mathcal{L}^\ep , \varphi \rangle | = |\int_{\Oe } \{ \nabla \mathbf{v}^\ep
- \, H(x_2)  \nabla \mathbf{u}^{eff}  +\notag \\  H(-x_2)
\nabla \sum_{j=1}^2 \mathbf{w}^{j, \ep} \frac{\p P^D }{ \p x_j}  \} \nabla \varphi \ dx -
 \int_{\Oe} \{  p^{\ep }
- H(x_2)  p^{eff} -\notag \\
 H(-x_2)
(\ep^{-2} P^D -\, \ep^{-1} \sum_{j=1}^2 \pi^{j, \ep} \frac{\p P^D}{ \p x_j} )
\} \; \mbox{ div }  \varphi  \ dx + \int_{ \{ x_2 =-H \} } \sum_{j=1}^2 \frac{\p }{\p x_2} (w_1^{j, \ep} \frac{\p P^D}{ \p x_j}) \varphi_1 \ dS \notag \\
 -\int_{\Sigma}  \bigl(  \sigma_0 \varphi
+ \sum_{j=1}^2 ( \nabla \mathbf{w}^{j, \ep}- \ep^{-1} \pi^{j, \ep}
I ) \frac{\p P^D}{ \p x_j} )   \mathbf{e}^2 \varphi  \bigr) \ dS |
\leq C \Vert \nabla \varphi \Vert _{L^2 (\Oe_2)^4 }  , \label{2.3A}
\end{gather}
where $\sigma_0 = (\nabla \mathbf{u}^{eff} -p^{eff} I )\mathbf{e}^2$ on $\Sigma$. 
\end{proposition}
{\bf Proof of proposition  \ref{Prop1}} We start with the weak formulation corresponding to
(\ref{1.3})-(\ref{1.8}):
\begin{equation} \label{2.2} \int_{\Oe } \nabla \mathbf{v}^\ep \nabla \varphi
 -\int_{\Oe }  p^{\ep } \; \mbox{ div }  \varphi =0, \quad \forall \varphi \in
V_{per} (\Oe )  .\end{equation}
As a first step we eliminate the boundary conditions.
The weak formulation corresponding to system (\ref{4.91})-(\ref{4.95}) is
\begin{equation}\label{W491}
  \int_{\Omega_1 } \nabla \mathbf{u}^{eff} \nabla \varphi
 -\int_{\Omega_1 }  p^{eff } \; \mbox{ div }  \varphi =-\int_\Sigma \sigma_0 \varphi \ dS, \quad \forall \varphi \in
V_{per} (\Oe )  .
\end{equation}
Next the weak formulation corresponding to the correction in the pore space $\Oe_2$ is
 \begin{gather} \label{WPM} \int_{\Oe_2 }  (- \nabla \sum_{j=1}^2 \mathbf{w}^{j, \ep} \frac{\p P^D }{ \p x_j} - \ep^{-2} P^D I +\, \ep^{-1} \sum_{j=1}^2 \pi^{j, \ep} \frac{\p P^D}{ \p x_j} I) \nabla \varphi \ dx
  =\int_{\Oe_2} (\mathbf{w}^{j, \ep} \Delta \frac{\p P^D }{ \p
x_j} +\notag \\
 2 \nabla \mathbf{w}^{j, \ep} \nabla \frac{\p P^D }{ \p x_j} -
\ep^{-1 }\pi^{j, \ep} \nabla \frac{\p P^D}{ \p x_j}) \cdot \varphi \ dx  - \int_\Sigma \sum_{j=1}^2 \bigg\{ ( \nabla \mathbf{w}^{j, \ep}- \ep^{-1} \pi^{j, \ep}
I ) \frac{\p P^D}{ \p x_j}    \mathbf{e}^2 +\notag \\
 \frac{\p^2 P^D}{\p x_2 \p x_j} \mathbf{w}^{j, \ep} \bigg\} \cdot \varphi   \ dS  +\int_{ \{ x_2 =-H \} } \sum_{j=1}^2 \frac{\p }{\p x_2} (w_1^{j, \ep} \frac{\p P^D}{ \p x_j}) \varphi_1 \ dS , \quad \forall \varphi \in
V_{per} (\Oe )  .\end{gather}
We observe that difference between (\ref{2.2}) and (\ref{W491})-(\ref{WPM}) is equivalent to
\begin{gather}
\langle \mathcal{L}^\ep , \varphi \rangle
= \int_{\Oe_2} \varphi \sum_{j=1}^2  A^j_\ep +\int_\Sigma \sum_{j=1}^2 \mathbf{w}^{j, \ep} \frac{\p^2 P^D}{\p x_2 \p x_j} \cdot \varphi \ dS, \label{2.3}\end{gather}
where  quantities $A^j_\ep $ are given by
\begin{gather}
A^j_\ep =  - \mathbf{w}^{j, \ep} \Delta \frac{\p P^D }{ \p
x_j} - 2 \nabla \mathbf{w}^{j, \ep} \nabla \frac{\p P^D }{ \p x_j} +
\ep^{-1 }\pi^{j, \ep} \nabla \frac{\p P^D}{ \p x_j}  = \notag \\
  \mathbf{w}^{j, \ep} \Delta \frac{\p P^D }{ \p
x_j}   +\ep^{-1 }\pi^{j, \ep} \nabla \frac{\p P^D}{ \p x_j}   -  2 \; \mbox{ div }  \ \{  \nabla
 \frac{\p P^D }{ \p x_j} \otimes\mathbf{w}^{j, \ep}  \} ,\, j=1,2. \label{2.4}\end{gather}
 We note that
 $$ -\nabla \sum_{j=1}^2 \mathbf{w}^{j} (\frac{x }{\ep}) \frac{\p P^D }{ \p x_j}  = -  \sum_{j=1}^2 \nabla \mathbf{w}^{j} (\frac{x }{\ep}) \frac{\p P^D }{ \p x_j} - \sum_{j=1}^2 \nabla
 \frac{\p P^D }{ \p x_j} \otimes\mathbf{w}^{j, \ep} . $$
and  a straightforward calculation yields
\begin{gather}
    |  \int_{\Oe_2} \varphi \sum_j A^j_\ep   | \leq C \Vert \nabla \varphi
\Vert _{L^2 (\Oe_2)^4 } \label{2.6} \\
| \int_\Sigma \sum_{j=1}^2 \mathbf{w}^{j} (\frac{x }{\ep}) \frac{\p^2 P^D}{\p x_2 \p x_j}  \cdot \varphi \ dS |\leq C \sqrt{\ep} \Vert \nabla \varphi
\Vert _{L^2 (\Oe_2)^4 } .\label{2.6A}
\end{gather}
\eop
\begin{remark}\label{Rno1} Now we see why it is necessary to impose $P^D =0$ at the interface $\Sigma$. Without it there would be a term $\displaystyle \int_\Sigma \ep^{-2} P^D \varphi_2 \ dS$ at the right hand side of (\ref{WPM}).
\end{remark}
\begin{remark} \label{Rno2} The candidate for the approximation of $\{ \mathbf{v}^\ep, p^{\ep} \}$ is
\begin{equation}\label{Approx0}
    \left\{
      \begin{array}{ll}
        \displaystyle \mathbf{v}^\ep \approx H(x_2)  \mathbf{u}^{eff}  - H(-x_2)
 \sum_{j=1}^2 \mathbf{w}^{j, \ep} \frac{\p P^D }{ \p x_j} ; &  \\
       \displaystyle  p^\ep \approx H(x_2)  p^{eff} +
 H(-x_2)
(\ep^{-2} P^D -\, \ep^{-1} \sum_{j=1}^2 \pi^{j, \ep} \frac{\p P^D}{ \p x_j} ). &
      \end{array}
    \right.
\end{equation}
Unfortunately, with such approximation we do not have continuity of the trace of the velocity approximation on the interface $\Sigma$.
\end{remark}
\subsection{Incremental accuracy correction, the 2nd part}\label{secpart}
The idea is to insert the correction to $\mathbf{v}^\ep$
as the test function $\varphi $ in equation (\ref{2.3}). Therefore the correction should be an element
of $ V_{per} (\Oe )$ and in this step we eliminate the trace jump on $\Sigma $. As in \cite{JaMi2}, fixing the traces on $\Sigma $ requires using the boundary layers defined by (\ref{BJ4.2})-(\ref{4.6}). At this stage we introduce the error functions
\begin{gather}
\mathbf{U}^\ep = \mathbf{v}^\ep - H(x_2)  (\mathbf{u}^{eff} - \mathbf{e}^1  C^{2, bl}_1 \frac{\p P^D }{ \p x_2} |_\Sigma ) + \notag \\
 H(-x_2)
 \sum_{j=1}^2 \mathbf{w}^{j, \ep} \frac{\p P^D }{ \p x_j}
-  \beta^{2,bl,\ep}  \frac{\p P^D }{ \p x_2} |_\Sigma ;\label{2.8}\\
P^\ep =   p^{\ep } - H(x_2) p^{eff}   - H(-x_2) (\ep^{-2} P^D
-\,  \ep^{-1} \sum_{j=1}^2 \pi^{j, \ep} \frac{\p P^D}{ \p x_j} ) \notag \\
-\ep^{-1}  ( \omega^{2,bl,\ep} -C^2_\pi ) \frac{\p P^D }{ \p x_2} |_\Sigma
, \label{2.9}\end{gather}
where
 $\{ \beta^{2, bl , \ep } , \omega^{2,
bl , \ep} \} (x) = \{ \beta^{2, bl } , \omega^{2,
bl } \} (\frac{x}{\ep})$ are defined by (\ref{BJ4.2})-(\ref{4.6}) and  $( C^{2, bl } , C^2_\pi )$ by (\ref{decay1}).
\begin{proposition}\label{Prop2}
 $\mathbf{U}^\ep \in H^1 (\Oe )^3$
and for all $\varphi \in V_{per} (\Oe ) $ we have
\begin{gather}
|\int_{\Oe} \nabla \mathbf{U}^\ep \nabla \varphi \ dx - \int_{\Oe} P^\ep  \mbox{ div }  \varphi \ dx + \int_{\{ x_2 =-H \} } \sum_{j=1}^2 \frac{\p }{\p x_2} (w_1^{j, \ep} \frac{\p P^D}{ \p x_j}) \varphi_1 \ dS| \leq \notag \\
C \Vert \nabla \varphi \Vert _{L^2 (\Oe_2)^4 }  + \Vert  \varphi \Vert _{H^1 (\Omega_1)^2 } .\label{Esttrac}
\end{gather}
\end{proposition}
{\bf Proof of proposition  \ref{Prop2} .} We have the following variational equation for $\{ \mathbf{U}^\ep , P^\ep \} $, for all $\varphi \in V_{per} (\Oe)$,:
\begin{gather}
\int_{\Oe} \nabla \mathbf{U}^\ep \nabla \varphi \ dx - \int_{\Oe} P^\ep \; \mbox{ div }  \varphi \ dx =
 \int_{\Sigma} \bigl(   \sigma_0  
+ \sum_{j=1}^2
  B^j_\ep   \bigr) \mathbf{e}^2 \varphi \ dS - \notag \\
 \int_{\{  x_2 =-H \} } \mathcal{C}^1_\ep \varphi_1 \ dS
+ \int_{\Oe_2} \varphi (\sum_{j=1}^2  A^j_\ep -  A^{22}_\ep ) \ dx
-2 \int_{\Oe }  A^{12}_\ep \nabla \varphi \ dx \notag \\
- \int_{\O_1}  (  A^{32}_\ep  + A^{42}_\ep )\varphi \ dx
, \label{2.10} \end{gather}
where
\begin{gather}
B^{j} = -\mathbf{w}^{j, \ep} \otimes \nabla \frac{\p P^D}{ \p x_j} , \label{2.11} \\
A^{12}_\ep = - \frac{d}{d x_1} (\frac{\p P^D}{ \p x_2} |_\Sigma ) \mathbf{e}^1 \otimes ( \beta^{2,bl,\ep} - H(x_2) C^{2, bl } )  ,  \label{2.13} \\
A^{22}_\ep = - \beta^{2,bl,\ep} \frac{d^2}{d x_1^2} (\frac{\p P^D}{ \p x_2} |_\Sigma )- \ep^{-1}
 ( \omega^{2,bl,\ep} -C^2_\pi ) \frac{d}{d x_1} (\frac{\p P^D}{ \p x_2} |_\Sigma ) \mathbf{e}^1 ,  \label{2.14} \\
A^{32}_\ep = - ( \beta^{2,bl,\ep} - C^{2, bl }_1 \mathbf{e}^1 ) \frac{d^2}{d x_1^2} (\frac{\p P^D}{ \p x_2}|_\Sigma ) , \label{2.15} \\
A^{42}_\ep = -\ep^{-1}  ( \omega^{2,bl,\ep} -C^2_\pi ) \frac{d}{d x_1} (\frac{\p P^D}{ \p x_2} |_\Sigma )\mathbf{e}^1 ,\label{2.16} \\
\mathcal{C}^1_\ep = \frac{\p {U}^\ep_1}{\p x_2} |_{\{ x_2=-H \} } = \sum_{j=1}^2 \frac{\p }{\p x_2} (w_1^{j, \ep} \frac{\p P^D}{ \p x_j})|_{\{ x_2=-H \} } +\, \mbox{exponentially small terms}.\label{2.16A}
\end{gather}
Then we have
\begin{gather}
\mid  \int_{\Sigma }  \sum_{j=1}^2 B^{j} \mathbf{e}^2 \varphi  \mid 
\leq C
\ep^{1/2} \Vert \nabla \varphi \Vert _{L^2 (\Oe_2)^4 } . \label{2.18}\end{gather}
Now we turn to the volume terms. We have
\begin{gather}
\mid  \int_{\Oe }  \sum_j A^{12}_\ep \nabla \varphi  \ dx \mid \leq C
\ep^{1/2} \Vert \nabla \varphi \Vert _{L^2 (\Oe)^4 } \label{2.19} \\
\mid  \int_{\Oe_2 }  \sum_j A^{22}_\ep \varphi  \ dx \mid \leq C
\Vert \nabla \varphi \Vert _{L^2 (\Oe_2)^4 } \label{2.20} \\
\mid  \int_{\O_1 }  \sum_j A^{32}_\ep \varphi  \ dx \mid \leq C\sqrt{\ep} 
\Vert \varphi \Vert _{L^2 (\O_1)^2 } . \label{2.21}\end{gather}
Finally, we estimate the term involving $A^{42}_\ep $.
Let $Q^{2}$ be defined by
\begin{equation}\label{(3.59P)}
    \left\{
      \begin{array}{ll}
       \displaystyle \frac{\p Q^2 }{ \p y_1 } = \omega ^{2,bl} - C^2_\pi  , & \hbox{on } (0,1)\times (0, +\infty ) ;\\
       \displaystyle Q^2   & \hbox{is } y_1-\hbox{periodic}.
      \end{array}
    \right.
\end{equation}
By definition of $ C^2_\pi$,  the function
\begin{equation}\label{(3.60P)}
 Q^2 (y_1 , y_2 ) =\int_0^{y_1} \omega^{2, bl} (t, y_2) dt - C^2_\pi y_1,
\quad y\in (0,1)\times (0, +\infty )
\end{equation}
is a solution for (\ref{(3.59P)}) and, using the results from \cite{JaMi2}, page 459,  there exists a constant $\gamma_0 >0$ such that
$\displaystyle e^{\gamma_0  y_2 }  Q^j \in L^2 (Z^+ ) .$

We set $Q^{2,\ep} (x) = \ep Q^{2} (x/\ep) $, $x\in \Omega_1$.
 Then  we obtain
\begin{equation}\label{Qprop}
  \left\{
    \begin{array}{ll}
      \displaystyle \frac{\p Q^{2,\ep }}{\p x_1} = \omega^{ 2, bl, \ep} (x) - C^2_\pi ;  &  \\
      \displaystyle  \Vert Q^{2,\ep } \Vert_{L^2 (\Omega_1)} \leq C\ep^{3/2} . &
    \end{array}
  \right.
\end{equation}
Therefore we have
\begin{gather}
 \mid \int_{\O_1}   A^{42}_\ep  \varphi \ dx \mid =  \mid \int_{\O_1}
\ep^{-1} Q^{2,\ep} (\varphi_1 \frac{d^2  }{ d x^2_1 } (\frac{\p P^D}{ \p x_2} |_\Sigma ) + \frac{\p  \varphi_1 }{ \p x_1 } \frac{\p }{\p x_1}  (\frac{\p P^D}{ \p x_2} |_\Sigma )) \mid \notag \\
\leq C  \ep^{1/2} \Vert  \varphi \Vert _{H^1 (\O_1)^2 }
 . \label{2.22}\end{gather}
Now
  the estimates
(\ref{2.18}) - (\ref{2.22}) show that the right hand side
in (\ref{2.10}) is bounded by
$$ C  \Vert \nabla \varphi \Vert _{L^2 (\Oe_2)^4 }
+ C \ep^{1/2} \Vert  \varphi \Vert _{H^1 (\O_1)^2 }
. $$\eop
\begin{remark}\label{Rno3} We would like to use $\mathbf{U}^\ep $ as a test function in the variational equation (\ref{2.10}). The difficulty with $\mathbf{U}^\ep $ is that the boundary condition at $\{ x_2 =-H \}$ is not satisfied. Hence we have to adjust its values at that boundary.
\end{remark}

\subsection{Incremental accuracy correction, the 3rd part: correction of the outer boundary effects}

First we calculate values of ${U}^\ep_2$ and $\displaystyle \frac{\p }{\p x_2} U_1^\ep$ at the lower outer boundary $\{ x_2 =-H \} $. We have
\begin{gather*}
  {U}^\ep_2 (x_1 , -H ) = v_2^\ep (x_1 , -H ) + \sum_{j=1}^2  \frac{\p P^D}{ \p x_j} (x_1 , -H ) K_{2j} + \sum_{j=1}^2  \frac{\p P^D}{ \p x_j} (x_1 , -H ) (w_2^j (\frac{x_1}{\ep} , 0) - K_{2j} ) \\
 + O(e^{Cx_2 /\ep } ) =
 \sum_{j=1}^2  \frac{\p P^D}{ \p x_j} (x_1 , -H ) (w_2^j (\frac{x_1}{\ep} , 0) - K_{2j} ) + \, \mbox{exponentially small terms}, \\
\frac{\p }{\p x_2} U_1^\ep (x_1 , -H ) = \frac{1}{\ep} \sum_{j=1}^2  \frac{\p P^D}{ \p x_j} (x_1 , -H ) \frac{\p w^j_1}{\p y_2} (\frac{x_1}{\ep} , 0) + \sum_{j=1}^2  \frac{\p^2 P^D}{ \p x_j \p x_2} (x_1 , -H ) w_1^j (\frac{x_1}{\ep} , 0) \\
+ \, \mbox{exponentially small terms}.
\end{gather*}

We follow again \cite{JaMi2} and correct the outer boundary effects using the corresponding boundary layer:
\begin{gather}
-\triangle  \mathbf{q}^{j,bl} +\nabla z^{j,bl}= 0 \qquad \hbox{ in }  Z^- \label{(3.51O)}\\
\hbox{div}   \ \mathbf{q}^{j ,bl} =0\qquad \hbox{ in }  Z^- \label{(3.52O)}\\
 {q}_2^{j,bl} =  K_{2j} - w_2^j \quad \mbox{and} \quad \frac{\p q_1^j}{\p y_2} = -\frac{\p w_1^j}{\p y_2} \hbox{ on } S \label{(3.53O)}\\
 \mathbf{q}^{j,bl} =0 \quad \hbox{ on } \displaystyle\cup_{k=1}^\infty \{ \p Y_F \setminus \p Y- (0,k) \},
 \quad  \{ \mathbf{q}^{j,bl} , z^{j,bl} \} \, \hbox{ is } y_1-\hbox{periodic.}
\label{(3.55O)}\end{gather}
Following the theory from \cite{JaMi2},  problem (\ref{(3.51O)})-(\ref{(3.55O)}) admits a unique solution $\mathbf{q}^{j,bl} \in H^1 (Z^-)^2$, smooth in $Z^-$. Furthermore, there is $\gamma_0 >0$ such that
$e^{\gamma_0 | y_2 | } \mathbf{q}^{j, bl} \in L^2 (Z^-)^2$ and, after adjusting a constant,  $e^{\gamma_0 | y_2 | } z^{j, bl} \in L^2 (Z^-)$.
The new error functions read
\begin{gather}
\mathbf{U}^{1, \ep } = \mathbf{U}^{ \ep }  + \sum_{j=1}^2  \frac{\p P^D}{ \p x_j} (x_1 , -H ) \mathbf{q}^{j, bl} (\frac{x_1}{\ep} , -\frac{x_2 +H}{\ep} ) , \label{Errorout1} \\
  P^{1, \ep } = {P}^{\ep }  + \frac{1}{\ep} \sum_{j=1}^2  \frac{\p P^D}{ \p x_j} (x_1 , -H ) z^{j, bl} (\frac{x_1}{\ep} , -\frac{x_2 +H}{\ep} ) . \label{Errorout2}
\end{gather}
Variational equation (\ref{2.10}) becomes
\begin{gather}
\int_{\Oe} \nabla \mathbf{U}^{1,\ep } \nabla \varphi \ dx - \int_{\Oe} P^{1, \ep} \; \mbox{ div }  \varphi \ dx =
 \int_{\Sigma} \bigl(   \sigma_0   + \sum_{j=1}^2
  B^j_\ep
   \bigr) \mathbf{e}^2 \varphi \ dS - \notag \\
 \int_{\{  x_2 =-H \} } \sum_{j=1}^2 w^{j, \ep}_1  \varphi_1 \frac{\p^2 P^D}{ \p x_j \p x_2}  \ dS
+ \int_{\Oe_2} \varphi (\sum_{j=1}^2  A^j_\ep -  A^{22}_\ep ) \ dx
-2 \int_{\Oe }   A^{12}_\ep \nabla \varphi \ dx \notag \\
- \int_{\O_1}  (  A^{32}_\ep  + A^{42}_\ep )\varphi \ dx -2 \int_{\Oe_2} \sum_{j=1}^2 \nabla \mathbf{q}^{j, bl} (\frac{x_1}{\ep} , -\frac{x_2 +H}{\ep} ) \nabla \frac{\p P^D}{ \p x_j} (x_1 , -H ) \varphi \ dx - \notag \\   \int_{\Oe_2} \sum_{j=1}^2 ( \Delta \frac{\p P^D}{ \p x_j} (x_1 , -H ) \mathbf{q}^{j, bl} (\frac{x_1}{\ep} , -\frac{x_2 +H}{\ep} ) +\notag  \\
\frac{1}{\ep} \nabla \frac{\p P^D}{ \p x_j} (x_1 , -H ) z^{j, bl} (\frac{x_1}{\ep} , -\frac{x_2 +H}{\ep} )) \varphi \ dx.
\label{2.10OUT} \end{gather}
The form of the right hand side of variational equation (\ref{2.10OUT}) yields
\begin{proposition}\label{Prop3}
We have $\mathbf{U}^{1,\ep } \in V_{per} (\Oe )$
and $\forall \varphi \in V_{per} (\Oe ) $ we have
\begin{gather}
|\int_{\Oe} \nabla \mathbf{U}^{1,\ep} \nabla \varphi \ dx - \int_{\Oe} P^{1,\ep}  \mbox{ div }  \varphi \ dx | \leq
C (\Vert \nabla \varphi \Vert _{L^2 (\Oe_2)^4 }  + \Vert  \varphi \Vert _{H^1 (\Omega_1)^2 } ).\label{EsttracBdr}
\end{gather}
\end{proposition}

\begin{remark}\label{Rno4}
It remains to estimate the pressure through the velocity and then to use the velocity error as a test function in equation (\ref{2.10}).
However at this stage the difficulties are coming from the compressibility effects in the term $\displaystyle
 \int_{\Oe} P^\ep  \mbox{ div }  \mathbf{U}^\ep \ dx$. In fact
\begin{gather*}
\; \mbox{ div}  \ \mathbf{U}^{1,\ep }= H(-x_2) \sum_{j=1}^2 \mathbf{w}^{j, \ep} \nabla \frac{\p P^D}{ \p x_j}
+( \beta^{2,bl,\ep}_1 - H(x_2) C^{2,bl}_1 ) \frac{d }{d x_1} (\frac{\p P^D}{ \p x_2} |_\Sigma ) \\
+  \sum_{j=1}^2  \frac{d}{d x_1} \frac{\p P^D}{ \p x_j} (x_1 , -H ) {q}_1^{j, bl} (\frac{x_1}{\ep} , -\frac{x_2 +H}{\ep}  )  \end{gather*}
and the estimate of the divergence is
 $\displaystyle \Vert  \mbox{ div }  \mathbf{U}^{1, \ep }\Vert _{L^2 (\Oe ) } \leq C.$
Therefore, we have to diminish the value of  $ \; \mbox{ div }  \ \mathbf{U}^{1, \ep} $.
\end{remark}
\subsection{Incremental accuracy correction, the 4th step: correction of the compressibility effects}

We start by introducing the correction the basic auxiliary problem, linked to the permeability auxiliary problems:
\begin{equation}\label{1.18}
    \left\{
      \begin{array}{ll}
        \displaystyle  \mbox{div}_y \gamma^{j,i} = w^j_i - \frac{K_{ij} }{ \mid Y_F \mid } \quad \mbox{in} \quad Y_F ; &  \\
        \displaystyle \gamma^{j,i} =0 \quad \hbox{ on } \p Y_F \setminus \p Y, \quad  \gamma^{j,i}
\quad \mbox{is} \quad 1-\mbox{periodic.}
 &
      \end{array}
    \right.
\end{equation}

The existence of at least one $\gamma^{j,i} \in H^1 ( Y_F)^2 \cap
C^{\infty}_{loc} \bigl( \cup_{k\in \mathbb{N}} (Y_F - (0,k) )^2 \bigr) $, satisfying
(\ref{1.18}) is straightforward.

We introduce $\gamma^{j,i,\ep }$ by
\begin{equation}\label{1.19}
  \gamma^{j,i, \ep } (x)= \ep \gamma^{j,i} (x/ \ep ), \quad x\in \Oe_2
\end{equation}
and extend it by zero to $\O_2 \setminus \Oe_2 .$

$\gamma^{j,i, \ep } $ is defined only in the porous part $\Omega_2$ and an auxiliary boundary layer velocity and pressures, correcting its values of on $\Sigma $, is needed.

First we construct  $\{ \gamma^{j,i,bl} , \pi ^{j,i,bl} \} $ satisfying
\begin{gather}
-\triangle _y \gamma^{j,i,bl} +\nabla_y \pi ^{j,i,bl}=0\qquad \hbox{ in } Z^+ \cup
Z^- , \label{1.40} \\
\mbox{div}_y  \gamma^{j,i,bl} =0\qquad \hbox{ in } Z^+ \cup Z^- , \label{1.41}\\
\bigl[ \gamma^{j,i,bl} \bigr]_S (\cdot , 0) = \gamma^{j,i}
(\cdot , 0) \quad \hbox{ on } S,  \label{1.42}\\
\bigr[ \{ \nabla_y \gamma^{j,i,bl} -\pi^{j,i,bl} I \} \mathbf{e}^2 \bigl]_S (\cdot , 0)
= \nabla_y \gamma^{j,i} (\cdot , 0) \mathbf{e}^2 \ \hbox{ on } S, \label{1.43} \\
\gamma^{j,i,bl} =0 \quad \hbox{ on } \displaystyle\cup_{k=1}^\inf \{ \p Y_F \setminus \p Y
-(0,k) \}, \quad  \{ \gamma^{j,i,bl} , \pi^{j,i,bl} \} \, \hbox{ is } y_1-
\hbox{periodic.} \label{1.44} \end{gather}
Proposition 3.19    from \cite{JaMi2} gives the existence of a
solution  $\{ \gamma^{j,i,bl}, \pi^{j,i,bl} \} \in V\cap
C^{\inf}_{\hbox{loc} } (Z^+ \cup Z^-)^2 \times C^{\inf}_{\hbox{loc} } (Z^+
\cup Z^-)  $ to equations (\ref{1.40})-(\ref{1.44}), where $\gamma^{j,i,bl} $ is uniquely determined
and $ \pi^{j,i,bl} $ is unique up to a constant.
$\displaystyle   \gamma^{j,i,bl} (\cdot , \pm 0 ) \in W^{2 - 1/ q , q } (S)^2 $ and
$\displaystyle   \{ \nabla \gamma^{j,i,bl} -\pi^{j,i,bl} I \} (\cdot , \pm 0 ) e_2
\in W^{1 - 1/ q , q } (S)^2 $, $\forall q\in [1, \infty [ $,
but the limits from two sides of $S$ are in general different.
Furthermore, it is proved that there exist  constants $\gamma_0 \in ]0,1[ ,
C^{j,i}_\pi $, 
and a constant vector $C^{j,i,bl}$ such that
$$ e^{\gamma_0 \mid y_2 \mid} \nabla_y \gamma^{j,i,bl} \in L^2 (Z^+ \cup Z^- )^4, \,
e^{\gamma_0 \mid y_2 \mid} \gamma^{j,i,bl} \in L^2 ( Z^- )^2  ,
e^{\gamma_0 \mid y_2 \mid} (\pi^{j,i,bl} - C^{j,i}_\pi )\in L^2 ( Z^+ )
$$
and
\begin{equation}\label{1.45}
    \left\{
      \begin{array}{ll}
        \displaystyle \mid \gamma^{j,i,bl} (y_1 , y_2 ) - \mathbf{C}^{j,i,bl} \mid  \leq
Ce^{-\gamma_0  y_2 } , \qquad y_2 > y_* ; & \\
        \displaystyle \mid \pi^{j,i,bl} ( y_1 , y_2 ) - C^{j,i}_\pi \mid  \leq
Ce^{-\gamma_0  y_2 } , \qquad y_2 > y_*  .&
      \end{array}
    \right.
\end{equation}
We define
\begin{equation}\label{1.46}
    \gamma^{j,i, bl, \ep} (x)= \ep \gamma^{j,i, bl} ({x\over \ep } ) \qquad
\hbox{and} \qquad \pi^{j,i,bl,\ep} (x)= \pi^{j,i,bl} ({x\over \ep } ),
\quad x\in \Oe,
\end{equation}
and extend $\gamma^{j,i, bl ,\ep } $ by zero to  $\O \setminus \Oe $.

Next, we need a correction for the compressibility effects coming from the boundary layer term $\displaystyle
 \beta^{2,bl,\ep}  \frac{\p P^D }{ \p x_2} |_\Sigma $:
We look for $\theta^{bl} $ satisfying
\begin{equation}\label{3.68}
    \left\{
      \begin{array}{ll}
        \displaystyle \mbox{div }  \theta^{ bl} = \beta^{2, bl }_1 -  C^{2, bl}_1 H(y_2)
 \quad  \mbox{in} \;   Z^+ \cup Z^-  ; &  \\
        \displaystyle \bigl[ \theta^{ bl} \bigr]_S = (\int_{Z_{BL}} ( C^{2, bl}_1 H(y_2) -
 \beta^{2, bl }_1 ) \ dy ) \mathbf{e}^2 \quad \mbox{on} \; S; & \\
        \displaystyle \theta^{bl}  = 0 \quad \hbox{ on } {\displaystyle\cup_{k=1}^\inf \{ \p
Y_F \setminus \p Y  -(0,k) \} }, \quad \theta^{bl} \quad \mbox{is} \;  y_1-\mbox{periodic.} &
      \end{array}
    \right.
\end{equation}
 After proposition 3.20 from \cite{JaMi2},
problem (\ref{3.68}) has at least one solution $\theta^{bl} \in H^1
(Z^+ \cup Z^- )^2 \cap C^{\infty }_{loc} (Z^+ \cup Z^- )^2.$ Furthermore,
$\theta^{bl} \in W^{1, q} ((0,1)^2)^2 $ and $ \theta^{bl} \in
W^{1, q} (Y- (0,1) )^2 $,  $\forall q\in [1, \infty ).$ Furthermore, there is $\gamma_0 >0$ such that
$\displaystyle e^{\gamma_0 \mid y_2 \mid } \theta^{j,i,bl} \in H^1 (Z^+ \cup Z^- )^2 .$

Let $\gamma^{j,i, \ep }$ be defined by (\ref{1.19})  and  $\gamma^{j,i, bl, \ep } ,
\pi^{j,i, bl, \ep } $, $ C^{j,i}_\pi , \mathbf{C}^{j,i, bl} $  by (\ref{1.45})-(\ref{1.46}).  We modify
$\{ \mathbf{u}^{eff} , p^{eff} \} $ by adding to it $\ep \{ \mathbf{u}^{1, eff} , p^{1, eff} \} $, satisfying (\ref{4.91})-(\ref{4.94}) and  with (\ref{4.95}) replaced by
\begin{gather}
    u^{1, eff}_2 =  - \ \sum_{j,k=1}^2 C^{j,k, bl}_2 \frac{\p^2 P^D}{\p x_j \p x_k} |_\Sigma  -  \theta^{ bl}_2 (\frac{x}{\ep}) \frac{d}{d x_1} \frac{\p P^D}{\p x_2} |_\Sigma  \quad \hbox{ on } \quad \Sigma  , \label{4.95AA} \\
   u^{1, eff}_1 = -  \sum_{j,k=1}^2 C^{j,k, bl}_1 \frac{\p^2 P^D}{\p x_j \p x_k} |_\Sigma -  \theta^{ bl}_1 (\frac{x}{\ep}) \frac{d}{d x_1} \frac{\p P^D}{\p x_2} |_\Sigma   \quad \hbox{ on } \quad \Sigma  . \label{4.95A}\end{gather}
The pair $\{ \mathbf{u}^{1, eff} , p^{1, eff} \} $ is uniquely defined by (\ref{4.91})-(\ref{4.94}), (\ref{4.95AA})-(\ref{4.95A}).

Finally we correct the compressibility effects coming from the boundary layer around $\{ x_2 = -H \} .$ We introduce
$\mathbf{Z}^{j, bl} $, $j=1,2$, satisfying
\begin{equation}\label{3.68C}
    \left\{
      \begin{array}{ll}
        \displaystyle \mbox{div}_y \  \mathbf{Z}^{j, bl} = q^{j, bl}_1
 \quad  \mbox{in} \;    Z^-  ; &  \\
        \displaystyle \bigl[ \mathbf{Z}^{j, bl} \bigr]_S = -(\int_{Z^{-}}
 q^{j, bl }_1 ) \ dy ) \mathbf{e}^2 \quad \mbox{on} \; S; & \\
        \displaystyle  \mathbf{Z}^{j, bl} = 0 \quad \hbox{ on } {\displaystyle\cup_{k=1}^\inf \{ \p
Y_F \setminus \p Y  -(0,k) \} }, \quad \mathbf{Z}^{j, bl} \quad \mbox{is} \;  y_1-\mbox{periodic.} &
      \end{array}
    \right.
\end{equation}
 After proposition 3.20 from \cite{JaMi2},
problem (\ref{3.68}) has at least one solution $\mathbf{Z}^{j, bl} \in H^1
(Z^+ \cup Z^- )^2 \cap C^{\infty }_{loc} (Z^+ \cup Z^- )^2.$ Furthermore,
$\mathbf{Z}^{j, bl} \in W^{1, q} ((0,1)^2)^2 $ and $ \mathbf{Z}^{j, bl} \in
W^{1, q} (Y- (0,1) )^2 $,  $\forall q\in [1, \infty ).$ Furthermore, there is $\gamma_0 >0$ such that
$\displaystyle e^{\gamma_0 \mid y_2 \mid } \mathbf{Z}^{j, bl} \in H^1 (Z^+ \cup Z^- )^2 .$
 Note that $\int_{Z^{-}}
 q^{j, bl }_2  \ dy =0$, $j=1,2$. Next we set
 $$ \mathbf{Z}^{j, bl, \ep} (x) = \ep \mathbf{Z}^{j, bl} (\frac{x_1}{\ep}, - \frac{x_2 +H}{\ep} ) , \quad x\in \Oe .$$

Now we introduce new velocity-pressure error functions by
\begin{gather}
\mathbf{U}^{2,\ep}  =\mathbf{U}^{1,\ep } - H(-x_2) \sum_{i,j=1}^2 \gamma^{j,i, \ep } \frac{\p^2 P^D}{ \p x_i \p x_j } - \sum_{i,j=1}^2 \bigl(  \gamma^{j,i, bl, \ep} -\ep  \mathbf{C}^{j,i, bl} H(x_2) \bigr) \frac{\p^2 P^D}{ \p x_i \p x_j } |_\Sigma  \notag \\
-H(x_2) \ep \mathbf{u}^{1, eff}  -\ep \theta^{bl}  (\frac{x}{\ep}) \frac{d}{d x_1} \frac{\p P^D}{\p x_2} |_\Sigma -\sum_{j=1}^2   \frac{d}{d x_1} \frac{\p P^D}{ \p x_j} (x_1 , -H ) ( \mathbf{Z}^{j, bl, \ep}  +\notag \\
\ep (\int_{Z^{-}}
 q^{j, bl }_1 ) \ dy ) R_\ep (\mathbf{e}^2 ) ) ;\label{ErrV2} \\
 P^{2, \ep} = P^{1, \ep} -H(x_2) \ep p^{1, eff} - \sum_{i,j=1}^2 \bigl( \pi^{j,i, bl, \ep } -  C^{j,i}_\pi \bigr) \frac{\p^2 P^D}{ \p x_i \p x_j } |_\Sigma
 , \label{ErrP2} \end{gather}
where $R_\ep$ is Tartar's restriction operator (see \cite{Ta1980}), defined after (\ref{Estdiv}).
\begin{proposition}\label{Prop4}
We have $\mathbf{U}^{2,\ep } \in V_{per} (\Oe )$
and for all $\varphi \in V_{per} (\Oe ) $ we have
\begin{gather}
|\int_{\Oe} \nabla \mathbf{U}^{2,\ep} \nabla \varphi \ dx - \int_{\Oe} P^{2,\ep}  \mbox{ div }  \varphi \ dx | \leq
C \Vert \nabla \varphi \Vert _{L^2 (\Oe_2)^4 }  + \Vert  \varphi \Vert _{H^1 (\Omega_1)^2 } ,\label{EsttracBdr22} \\
\Vert \mbox{div} \, \mathbf{U}^{2, \ep} \Vert_{L^2 (\Oe)} \leq C\ep .\label{divEst}
\end{gather}
\end{proposition}
{\bf Proof of proposition \ref{Prop4} :}
We prove that $\mathbf{U}^{2,\ep } \in V_{per} (\Oe )$ by a direct verification. Furthermore,
\begin{gather}
\mbox{div} \, \mathbf{U}^{2, \ep} = - H(-x_2) \sum_{i,j=1}^2
\gamma^{j,i, \ep } \nabla \frac{\p^2 P^D }{ \p x_i \p x_j}  -
 \sum_{i,j=1}^2 \bigl( \gamma^{j,i, bl, \ep }_1 -\ep  C^{j,i, bl}_1 H(x_2) \bigr)
\frac{d}{d x_1} \frac{\p^2 P^D }{ \p x_i \p x_j} |_\Sigma \notag \\
 \hskip-10pt - \ep \theta^{bl}_1 (\frac{x}{\ep}) \frac{d^2}{d x_1^2} \frac{\p P^D }{ \p x_2} |_\Sigma
 -\sum_{j=1}^2   \frac{d^2}{d x_1^2} \frac{\p P^D}{ \p x_j} (x_1 , -H ) ( \mathbf{Z}^{j, bl, \ep}_1  +\ep (\int_{Z^{-}}
 q^{j, bl }_1 ) \ dy ) (R_\ep (\mathbf{e}^2 ))_1 )
  ,\label{div2} \end{gather}
which yields (\ref{divEst}).

It remains to estimate the
right hand side and prove (\ref{EsttracBdr22}):
\begin{gather}
\int_{\Oe} \nabla \mathbf{U}^{2, \ep} \nabla \varphi - \int_{\Oe} P^{2, \ep} \mbox{ div } \varphi =
 \int_{\Sigma} \bigl(   \sigma_0 + \ep \sigma_1  + \sum_{j=1}^2
  B^j_\ep  \bigr) \mathbf{e}^2 \varphi \ dS - \notag \\
 \int_{\{  x_2 =-H \} } \sum_{j=1}^2 w^{j, \ep}_1  \varphi_1 \frac{\p^2 P^D}{ \p x_j \p x_2}  \ dS
+ \int_{\Oe_2} \varphi (\sum_{j=1}^2  A^j_\ep -  A^{22}_\ep ) \ dx
-2 \int_{\Oe }   A^{12}_\ep \nabla \varphi \ dx \notag \\
- \int_{\O_1}  (  A^{32}_\ep  + A^{42}_\ep )\varphi \ dx -2 \int_{\Oe_2} \sum_{j=1}^2 \nabla \mathbf{q}^{j, bl} (\frac{x_1}{\ep} , \frac{x_2 +H}{\ep} ) \nabla \frac{\p P^D}{ \p x_j} (x_1 , -H ) \varphi \ dx - \notag \\   \int_{\Oe_2} \sum_{j=1}^2 ( \Delta \frac{\p P^D}{ \p x_j} (x_1 , -H ) \mathbf{q}^{j, bl} (\frac{x_1}{\ep} , \frac{x_2 +H}{\ep} ) +\notag  \\
\frac{1}{\ep} \nabla \frac{\p P^D}{ \p x_j} (x_1 , -H ) z^{j, bl} (\frac{x_1}{\ep} , \frac{x_2 +H}{\ep} )) \varphi \ dx
+\int_{\Oe_2}  \sum_{j,i=1}^2  A^{1,j,i}_\ep \nabla \varphi \ dx +\notag \\
 \int_{\Oe}
\sum_{j,i=1}^2 A^{2,j,i}_\ep  \varphi \ dx
+\sum_{j,i=1}^2 \int_{\Oe }  \{ A^{3,j,i}_\ep  +  A^{4,j,i}_\ep \} \nabla \varphi \ dx
-  \sum_{i,j=1}^2 \int_{\Sigma } \bigl( B^{1,j,i}_\ep  + B^{2,j,i}_\ep \bigr)
\mathbf{e}^2 \varphi  \ dS, \notag \\
- \int_{\Oe} \nabla ( \sum_{j=1}^2   \frac{d}{d x_1} \frac{\p P^D}{ \p x_j} (x_1 , -H ) ( \mathbf{Z}^{j, bl, \ep}  +
\ep (\int_{Z^{-}}
 q^{j, bl }_1 ) \ dy ) R_\ep (\mathbf{e}^2 ) )) \nabla \varphi \ dx +\notag \\
  \int_{\Oe} (\sum_{i,j=1}^2 \bigl( \pi^{j,i, bl, \ep } -  C^{j,i}_\pi \bigr) \frac{\p^2 P^D}{ \p x_i \p x_j } |_\Sigma ) \mbox{ div } \varphi \ dx, \label{2.25}\end{gather}
where
\begin{gather}
A^{1,j,i}_\ep = -\nabla \gamma^{j,i, \ep} \frac{\p^2 P^D }{ \p x_i
\p x_j}  - \gamma^{j,i, \ep} \otimes \nabla \frac{\p^2 P^D }{ \p x_i  \p x_j} ,  \label{2.26} \\
A^{2,j,i}_\ep = -(\gamma^{j,i, bl,\ep} - \ep H(x_2) \mathbf{C}^{j,i, bl} ) \frac{d^2}{d x_1^2}
 \frac{\p^2 P^D }{ \p x_i  \p x_j} |_\Sigma -\notag \\
 ( \pi^{j, i, bl,\ep} -C^{j,i}_\pi H(x_2) )  \frac{d}{d x_1} \frac{\p^2 P^D }{ \p x_i  \p x_j} |_\Sigma \mathbf{e}^1 ,\label{2.27} \\
A^{3,j,i}_\ep =  - 2  \{ ( \gamma^{j,i, bl, \ep} -\ep H(x_2) \mathbf{C}^{j,i, bl} ) \otimes
\nabla \frac{\p^2 P^D }{ \p x_i  \p x_j} |_\Sigma  \} ,\label{2.28} \\
A^{4,j,i}_\ep = - \ep \nabla \theta^{bl} (\frac{x}{ \ep }) \frac{d}{d x_1} \frac{\p P^D }{  \p x_2} |_\Sigma \mathbf{e}^1
  -  \ep \theta^{bl} (\frac{x}{ \ep }) \otimes \nabla \frac{\p^2 P^D }{ \p x_1  \p x_2} |_\Sigma
  \label{2.29} \end{gather} \begin{gather}
B^{1,j,i}_\ep = -\gamma^{j, i, \ep} (\cdot , -0) \otimes \nabla \frac{\p^2 P^D }{ \p x_i  \p x_j} |_\Sigma
 - \nabla \gamma^{j, i, \ep} |_\Sigma \frac{\p^2 P^D }{ \p x_i  \p x_j} |_\Sigma
  \label{2.30} \\
B^{2,j,i}_\ep = -\ep \mathbf{C}^{j,i, bl} \otimes \nabla \frac{\p^2 P^D }{ \p x_i  \p x_j} |_\Sigma  . \label{2.31}\end{gather}
Then
\begin{gather}
\mid \sum_{j,i=1}^2 \int_{\Oe_2} A^{1,j,i}_\ep \nabla \varphi \ dx \mid \leq C \ep^{1/2} \Vert
\nabla \varphi \Vert_{L^2 (\Oe_2 )^4} \label{2.32} \\
\mid \sum_{j,i=1}^2 \int_{\Oe_2} A^{2,j,i}_\ep \varphi \ dx \mid \leq C \ep^{3/2} \Vert \nabla
\varphi \Vert_{L^2 (\Oe_2 )^4} \label{2.33} \\
\mid \sum_{j,i=1}^2 \int_{\O_1} A^{2,j,i}_\ep \varphi \ dx \mid \leq C \ep^{1/2}
  \Vert  \varphi \Vert_{L^2 (\O_1)^2}
 \label{2.34} \\
\mid \sum_{j,i=1}^2 \int_{\Oe} A^{3,j,i}_\ep \nabla \varphi \ dx \mid \leq C \ep^{3/2}
\Vert \nabla \varphi \Vert_{L^2 (\Oe )^4} \label{2.35} \\
\mid \sum_{j,i=1}^2 \int_{\Oe} A^{4,j,i}_\ep \nabla \varphi \mid \leq C
\ep^{1/2}  \Vert \nabla \varphi \Vert_{L^2 (\Oe )^4} . \label{2.36}\end{gather}
and
\begin{gather}
\mid \sum_{j,i=1}^2 \int_{\Sigma} B^{1,j,i}_\ep \mathbf{e}^2 \varphi \ dS \mid \leq C
\ep^{1/2}\Vert \nabla \varphi \Vert_{L^2 (\Oe_2 )^4} \label{2.37} \\
\mid \sum_{j,i=1}^2 \int_{\Sigma} B^{2,j,i}_\ep \mathbf{e}^2 \varphi \ dS\mid \leq C
\ep^{3/2}\Vert \nabla \varphi \Vert_{L^2 (\Oe_2 )^4} . \label{2.38} \end{gather}

Proposition \ref{Prop4} is proved.\eop

\begin{corollary} \label{cor1}
We have
\begin{gather} \int_{\Oe} |\nabla \mathbf{U}^{2,\ep} |^2 \ dx \leq C\ep || P^{2,\ep} ||_{L^2 (\Oe)} +
C \Vert \nabla  \mathbf{U}^{2,\ep} \Vert _{L^2 (\Oe_2)^4 }  + \Vert  \mathbf{U}^{2,\ep} \Vert _{H^1 (\Omega_1)^2 } ,\label{EsttracBdrC}
\end{gather}
\end{corollary}

Hence at this point we need to estimate the pressure error $P^{2,\ep}$ using the velocity error $\mathbf{U}^{2,\ep}$.

\subsection{Pressure estimates}\label{Proofs2}
Following \cite{JaMi2} we consider the Stokes system
\begin{gather}
 -\Delta \mathbf{a}^\ep + \nabla \zeta^\ep  = \mathbf{M}^1_\ep + \mbox{ div } M_\ep^2   \quad  \hbox{ in }
         \hskip1pt  \Oe  ;
\label{5.1P} \\
       \hbox{div\  } \mathbf{a}^\ep =0  \quad \hbox{ in } \hskip3pt \Oe
       ; \label{5.11P} \\ 
 \mathbf{a}^\ep =0 \quad  \hbox{ on } \hskip3pt \p \Oe \setminus \p \O \quad \mbox{and on} \quad \{ x_2 =h \}  ; \label{5.12P} \\
        {a}_2^\ep =0 \quad  \hbox{and} \quad \frac{\p a^\ep_1}{\p x_2} = G^\ep  \quad  \hbox{ on } \hskip3pt  \{ x_2 =-H \} ; \label{5.13P} \\
       \{ \mathbf{a}^\ep , \zeta^\ep \}   \quad  \hbox{ is } L-\hbox{periodic in } \; x_1 ; \label{5.14P} \\
[ \mathbf{a}^\ep  ]_{\Sigma} =0 \quad  \hbox{and} \quad [(\nabla \mathbf{a}^\ep -\zeta^\ep I) \mathbf{e}^2 ]_{\Sigma} = G^\ep_\Sigma .\label{5.15P}
\end{gather}
We have
\begin{gather*}
    \int_{\Oe} \zeta^\ep  \ \hbox{div\  } \varphi \ dx =  \int_{\Oe} \nabla \mathbf{a}^\ep \nabla \varphi \ dx - \int_{\Oe} \mathbf{M}^1_\ep \varphi  \ dx + \int_{\Oe} M_\ep^2   \nabla \varphi \ dx  + \\
\int_{\Sigma} G^\ep_\Sigma  \varphi \ dS +   \int_{\{ x_2 = -H \} } G^\ep  \varphi \ dS, \qquad \forall \varphi \in V(\Oe),
\end{gather*}
which yields
\begin{gather}
    | \int_{\Oe} \zeta^\ep  \ \hbox{div\  } \varphi \ dx | \leq  C\bigg\{  \Vert \nabla \mathbf{a}^\ep \Vert_{L^2 (\Oe)^4} + ||  \mathbf{M}^1_\ep ||_{L^2 (\Omega_1)^2} + \ep ||  \mathbf{M}^1_\ep ||_{L^2 (\Oe_2)^2} + || M_\ep^2 ||_{L^2 (\Oe)^4}  \notag  \\
+ \sqrt{\ep} \bigg( || G^\ep_\Sigma ||_{L^2 (\Sigma) } +  || G^\ep ||  _{ L^2 (\{ x_2 = -H \}) } \bigg) \bigg\} || \nabla \varphi ||_{L^2 (\Oe)^4} . \label{Estdiv}
\end{gather}
At this point we need Tartar's restriction operator $R_\ep$ (see \cite{All97}, \cite{SP80} and \cite{Ta1980}). It is constructed for every pore on the following way:

Let $\gamma $ be a smooth curve, strictly contained in the cell $Y$, and enclosing the solid part $Y_s$. Let $Y_M$ be the domain between $\gamma$ and $\p Y_s$. Then using an intermediary nonhomogeneous Stokes system in $Y_M$ a linear  operator $R: H^1 (Y)^2 \to H^1 (Y_F)^2$ is constructed, such that
\begin{gather*}
     R \mathbf{u} (y) = \mathbf{u} (y) \quad \mathbf{for} \quad y\in Y\setminus ({\bar Y}_s \cup  Y_M ) , \quad R \mathbf{u} (y) = 0 \quad \mathbf{for} \quad y\in  Y_s , \\
\mathbf{u} = 0\; \mbox{ on }  Y_s \quad \Rightarrow R \mathbf{u}  = \mathbf{u} \quad \mbox{on} \quad Y, \quad
\mbox{div } \mathbf{u} = 0\; \mbox{ on }  Y \quad \Rightarrow \mbox{div } (R \mathbf{u})  = 0 \quad \mbox{on} \quad Y, \\
|| R \mathbf{u} ||_{H^1 (Y_F)^2} \leq C || \mathbf{u} ||_{H^1 (Y)^2} , \quad \forall \mathbf{u} \in H^1 (Y)^2.
\end{gather*}

Next the operator $R_\ep : H^1 (\Omega)^2 \to H^1 (\Oe)^2$ is defined by applying the operator $R$ to each $\ep (Y+ k)$ cell. After \cite{All97}, \cite{SP80} and \cite{Ta1980}, we have
\begin{gather}
     R_\ep \mathbf{u} (x) =  0 \quad \mathbf{for} \quad x\in   \Omega \setminus  \Oe , \quad
\mathbf{u} = 0\; \mbox{ on }   \Omega \setminus  \Oe  \; \Rightarrow R_\ep \mathbf{u}  = \mathbf{u} \quad \mbox{on} \quad \Oe, \notag \\
\mbox{div } \mathbf{u} = 0\; \mbox{ on }  \Omega \quad \Rightarrow \mbox{div } (R_\ep  \mathbf{u})  = 0 \quad \mbox{on} \quad \Oe, \notag \\
|| R_\ep \mathbf{u} ||_{L^2 (\Oe)^2} \leq C || \mathbf{u} ||_{L^2 (\Omega )^2 } + C\ep || \nabla \mathbf{u} ||_{L^2 (\Omega )^4 }, \quad \forall \mathbf{u} \in H^1 (\Omega)^2 , \label{RL2} \\
|| \nabla R_\ep \mathbf{u} ||_{L^2 (\Oe)^2} \leq \frac{C}{\ep} || \mathbf{u} ||_{L^2 (\Omega )^2 } + C || \nabla \mathbf{u} ||_{L^2 (\Omega )^4 }, \quad \forall \mathbf{u} \in H^1 (\Omega)^2 . \label{RLG2}
\end{gather}
Next we follow the calculation of Lipton and Avellaneda from \cite{LA} and extend the pressure to $\Omega_2$ by
\begin{equation}\label{Extpress}
    {\tilde \zeta}^\ep (x) = \left\{
                           \begin{array}{ll}
                             \displaystyle \zeta^\ep (x) \quad \mbox{ for } \quad x\in \Oe ; &  \\
                             \displaystyle \frac{1}{\ep^2 | Y_M |} \int_{\ep (Y_s - (k_1 , k_2 ))} \zeta^\ep (y) \ dy, \quad \mbox{ for } \quad x\in  \ep (Y_M - (k_1 , k_2 ) ).&
                           \end{array}
                         \right.
\end{equation}

The calculation of Lipton and Avellaneda gives
\begin{equation}\label{equpress}
    \int_{\Omega} {\tilde \zeta}^\ep  \ \hbox{div\  } \varphi \ dx = \int_{\Oe} \zeta^\ep  \ \hbox{div\  } (R_\ep \varphi ) \ dx, \quad \forall \varphi \in V(\Oe).
\end{equation}

\begin{proposition}\label{prop5} We have
\begin{gather}
    || {\tilde \zeta}^\ep  ||_{L^2 (\Omega )}   \leq  \frac{C}{\ep} \bigg\{  \Vert \nabla \mathbf{a}^\ep \Vert_{L^2 (\Oe)^4} + ||  \mathbf{M}^1_\ep ||_{L^2 (\Omega_1)^2} + \ep ||  \mathbf{M}^1_\ep ||_{L^2 (\Oe_2)^2} + || M_\ep^2 ||_{L^2 (\Oe)^4} \notag   \\
+ \sqrt{\ep} \bigg( || G^\ep_\Sigma ||_{L^2 (\Sigma) } +  || G^\ep || _{ L^2 (\{ x_2 = -H \} )} \bigg) \bigg\}  . \label{Estpressext}
\end{gather}
\end{proposition}
{\bf Proof of proposition \ref{prop5}: $\, $}
Let $g\in L^2 (\Omega)$. We set $h= g- \displaystyle \frac{1}{| \Omega|} \int_\Omega g \ dx$. Obviously $ \int_\Omega h \ dx =0$. Let
\begin{gather}
\hskip-15pt V_{per} (\Omega ) = \{ \mathbf{z}\in H^1 (\Omega )^2 :
\, 
\mathbf{z} =0 \quad \mbox{on} \quad \{ x_2 =h\}, \, 
 z_2 =0 \quad \mbox{on} \quad \{ x_2 =-H\} \notag \\
\hbox{ and } \mathbf{z} \, \hbox{ is L-periodic in } x_1 \hbox{ variable } \} .
\label{2.1A} \end{gather}
Then there exists  $\varphi \in  V_{per} (\Omega )$ such that div $\varphi =h$ in $\Omega$ and $\displaystyle || \varphi ||_{H^1 (\Omega)^2} \leq C || h ||_{L^2 (\Omega )}$, for all $h\in L^2_0 (\Omega)$.

Therefore we have
\begin{gather*}
    \int_\Omega {\tilde \zeta}^\ep h \ dx = \int_\Omega {\tilde \zeta}^\ep \ \mbox{div } \varphi \ dx =
\int_{\Oe} { \zeta}^\ep \ \mbox{div } (R_\ep \varphi )\ dx
\end{gather*}
and using (\ref{Estdiv}) and (\ref{RL2}), (\ref{RLG2}) yields
\begin{gather}
    | \int_{\Omega} {\tilde \zeta}^\ep  \ h \ dx | \leq  \frac{C}{\ep} \bigg\{  \Vert \nabla \mathbf{a}^\ep \Vert_{L^2 (\Oe)^4} + ||  \mathbf{M}^1_\ep ||_{L^2 (\Omega_1)^2} + \ep ||  \mathbf{M}^1_\ep ||_{L^2 (\Oe_2)^2} + || M_\ep^2 ||_{L^2 (\Oe)^4}  \notag  \\
+ \sqrt{\ep} \bigg( || G^\ep_\Sigma ||_{L^2 (\Sigma) } +  || G^\ep ||_{L^2 (\{ x_2 = -H \} )} \bigg) \bigg\} || \nabla \varphi ||_{L^2 (\Omega)^4} . \label{EstdivM}
\end{gather}
Since
$$ \int_\Omega ( {\tilde \zeta}^\ep -  \frac{1}{| \Omega|} \int_\Omega {\tilde \zeta}^\ep \ dy ) g \ dx = \int_\Omega ( {\tilde \zeta}^\ep -  \frac{1}{| \Omega|} \int_\Omega {\tilde \zeta}^\ep \ dy ) h \ dx =
\int_\Omega {\tilde \zeta}^\ep h \ dx,$$
we conclude that ${\tilde \zeta}^\ep - \displaystyle \frac{1}{| \Omega|} \int_\Omega {\tilde \zeta}^\ep \ dy $ satisfies  bound (\ref{Estpressext}).

For the mean we have
$$ 0 = \int_{\Omega_1}  {\tilde \zeta}^\ep \ dx = \int_{\Omega_1} ( {\tilde \zeta}^\ep -  \frac{1}{| \Omega|} \int_\Omega {\tilde \zeta}^\ep \ dy )  \ dx  + \int_\Omega  \frac{| \Omega_1 |}{| \Omega|} {\tilde \zeta}^\ep \ dx $$
and
\begin{gather}
|  \frac{1}{| \Omega|} \int_\Omega {\tilde \zeta}^\ep \ dx |\leq  \frac{1}{| \Omega|^{1/2}} || {\tilde \zeta}^\ep -  \frac{1}{| \Omega|} \int_\Omega {\tilde \zeta}^\ep \ dy ||_{L^2 (\Omega_1)^2}.
\label{Meanz}
\end{gather}
Estimate (\ref{Meanz}) implies bound (\ref{Estpressext}) for ${\tilde \zeta}^\ep$.
\eop

\subsection{Global energy estimate and proof of theorem \ref{Thbasic} }

{\bf Proof of theorem \ref{Thbasic}: $\ $}
Now  we choose $\varphi = \mathbf{U}^{2,\ep} $ as  test function in (\ref{2.25}).
Using 
estimates (\ref{2.32}) - (\ref{2.38})
and estimate (\ref{Estpressext}) from  proposition \ref{prop5}, we obtain
\begin{gather}
  \mid \int_{\Oe } \nabla \mathbf{U}^{2, \ep}  \nabla \mathbf{U}^{2, \ep} \mid \leq {C\over \ep }
\{ \Vert \nabla \mathbf{U}^{2, \ep} \Vert_{L^2 (\Oe_2 )^4} + C\} \Vert  \mbox{ div }  \ \mathbf{U}^{2,\ep}
\Vert_{L^2 (\Oe_2 )} +\notag \\
 C \Vert \nabla \mathbf{U}^{2, \ep} \Vert_{L^2 (\Oe_2 )^4} ,  \label{2.39}
\end{gather}
 which yields
\begin{gather}
|| \nabla \mathbf{U}^{2, \ep} ||_{L^2 (\Oe )^4} \leq C, \label{Erro1}
   \\
  || \mbox{ div } \mathbf{U}^{2, \ep} ||_{L^2 (\Oe )^4} +||  \mathbf{U}^{2, \ep} ||_{L^2 (\Oe_2 )^2} \leq C\ep \label{Erro2}
   \\
||  P^{2,\ep} ||_{L^2 (\Oe)} \leq \frac{C}{\ep} \label{Erro3}
\end{gather}
Hence estimates (\ref{Est3})-(\ref{Est4}) are proved.

It remains to prove estimates (\ref{Est1})-(\ref{Est2}).

First (\ref{Erro1})-(\ref{Erro2}) imply
\begin{equation}\label{Sigmaerr}
  ||  \mathbf{U}^{2, \ep} ||_{L^2 (\Sigma )^2} \leq C\sqrt{\ep }
\end{equation}
and estimate (\ref{Est2}) is proved.

Next we remark that in $\Omega_1$
 the error functions $ \mathbf{U}^{2, \ep} $ and $P^{2,\ep}$ satisfy the system
 \begin{equation}\label{5.1U}
    \left\{
      \begin{array}{l@{}}
      \displaystyle  -\Delta \mathbf{U}^{2, \ep} + \nabla P^{2, \ep}  = \mathbf{G}^{1, \ep} + \mbox{ div }  G^{2,\ep} \quad  \hbox{ in }
         \hskip1pt  \O_1  ;\\ \noalign{\vskip+4mm}
       \hbox{div\  } \mathbf{U}^{2, \ep} =\Lambda^\ep   \quad \hbox{ in } \hskip3pt \O_1
       ; \\ \noalign{\vskip+4mm}
       \mathbf{U}^{2, \ep} =\xi^\ep \quad  \hbox{ on } \hskip3pt  \Sigma; \quad \mathbf{U}^{2, \ep} =0 \quad  \hbox{ on } \hskip3pt  \{ x_2 =h \}  ; \\ \\
       \{ \mathbf{U}^{2, \ep} , P^{2, \ep} \}   \quad  \hbox{ is }L\hbox{-periodic in } \; x_1 ,
      \end{array}
    \right.
\end{equation}
where, after neglecting the boundary layer tails,
\begin{gather}
\Lambda^\ep =  -
 \sum_{i,j=1}^2 \bigl( \gamma^{j,i, bl, \ep }_1 -\ep  C^{j,i, bl}_1  \bigr)
\frac{d}{d x_1} \frac{\p^2 P^D }{ \p x_i \p x_j} |_\Sigma
 - \ep \theta^{bl}_1 (\frac{x}{\ep}) \frac{d^2}{d x_1^2} \frac{\p P^D }{ \p x_2} |_\Sigma
\label{Lambda} \\
\mathbf{G}^{1, \ep} = (\frac{Q^{2, \ep }}{\ep} \mathbf{e}^1 + \beta^{2,bl,\ep} - C^{2, bl }_1 \mathbf{e}^1   )\frac{d^2}{d x_1^2} (\frac{\p P^D}{ \p x_2} |_\Sigma )
 + \sum_{j,i=1}^2 A^{2,j,i}_\ep , \label{G1} \\
G^{2, \ep} = 
  \frac{d}{d x_1} (\frac{\p P^D}{ \p x_2} |_\Sigma )  \mathbf{e}^1 \otimes (2  \beta^{2,bl,\ep} - 2 C^{2, bl }_1 \mathbf{e}^1 - \frac{Q^{2, \ep }}{\ep}  \mathbf{e}^1 )
 + \sum_{j,i=1}^2 (A^{3,j,i}_\ep + A^{4,j,i}_\ep) .\label{G2}
\end{gather}
 The function $Q^{2, \ep }$  is given by (\ref{(3.60P)}) and,  for $i,j=1,2$, $A^{2,j,i}_\ep$ , $A^{3,j,i}_\ep$ and $A^{4,j,i}_\ep$ by (\ref{2.27})-(\ref{2.29}).

After (\ref{div2}), we have $ \displaystyle || \Lambda^\ep ||_{L^2 (\Omega_1 )} \leq C\ep^{3/2} $ . Using the basic theory of the Stokes system (see e.g. \cite{Tem})  there exists $\{ \mathbf{b} , \kappa \} \in H^1 (\Omega_1 )^2 \times L^2 (\Omega_1 )$, such that
\begin{equation}\label{5.1U1}
    \left\{
      \begin{array}{l@{}}
      \displaystyle  -\Delta \mathbf{b} + \nabla \kappa  = 0 \quad  \hbox{ in }
         \hskip1pt  \O_1  ;\\ \noalign{\vskip+4mm}
       \hbox{div\  } \mathbf{b} =\Lambda^\ep   \quad \hbox{ in } \hskip3pt \O_1
       ; \\ \noalign{\vskip+4mm}
       \mathbf{b}  \quad  \hbox{is given  on } \hskip3pt \Sigma_T = \Sigma \cup \{ x_2 =h \}  \quad \mbox{and}
       \quad  ||  \mathbf{b} ||_{H^{1/2} (\Sigma_T)^2} \leq C\ep^{3/2} ; \\ \\
       \{ \mathbf{b} , \kappa \}   \quad  \hbox{ is }L\hbox{-periodic in } \; x_1 ,
      \end{array}
    \right.
\end{equation}
Now we see that the pair $\{ \mathbf{U}^{2, \ep} - \mathbf{b}, P^{2, \ep} -\kappa \}$ satisfies system (\ref{5.1U}) with $\Lambda^\ep =0$. Such system admits the notion of a very weak solution, introduced by transposition. We refer to \cite{CO} , pages 61-68, and \cite{AMCAM2011} for the definition and properties of a very weak solution. Note that $\int_\Sigma (\xi^\ep_2 - b_2 ) \ dS=0$.

Let $H^k_p (\Omega_1 )^2 = \{ \mathbf{z} \in  H^k (\Omega_1 )^2 \ | \; \mathbf{z} \; \mbox{ is L-periodic} \; \mbox{ and } \; \mathbf{z} =0 \quad \mbox{on } \; \{ x_2 =h \} \  \}$, $k=1,2$.
Then the $q=r=2$-version of proposition 4.2., page 302, from \cite{AMCAM2011}, gives the estimate
 \begin{gather}
\hskip-19pt || \mathbf{U}^{2, \ep} - \mathbf{b} ||_{L^2 (\Omega_1 )^2 } \leq C \{ || \mathbf{G}^{1, \ep} ||_{(H^2_p (\Omega_1 )^2 )'} +
||{G}^{2, \ep} ||_{(H_p^1 (\Omega_1 )^4)' } + || \xi^\ep - \mathbf{b} ||_{L^2 (\Sigma_T)^2 } \}
\label{estauxom1}
\end{gather}
Using estimates (\ref{2.19}), (\ref{2.21})  and (\ref{2.34})-(\ref{2.36}), choosing $Q^{2, \ep}$ with zero mean
and repeating  calculations analogous to ones from (\ref{2.22}) to other terms, yields
\begin{equation}\label{auxforc}
  || \mathbf{G}^{1, \ep} ||_{(H^1_p (\Omega_1 )^2 )' } +
||{G}^{2, \ep} ||_{(H^1_p (\Omega_1 )^4)' } \leq C \ep^{3/2}.
\end{equation}
Now we are able to conclude that
\begin{gather}
|| \mathbf{U}^{2, \ep}  ||_{L^2 (\Omega_1 )^2 } \leq C \sqrt{\ep}
\label{estauxom2}\end{gather}
and estimate (\ref{Est1}) is proved.\eop

\section{Proof of theorem \ref{Thadvanced}}\label{Proofs3}

The proof is in fact a slight modification of the proof of theorem \ref{Thbasic}. Our goal is to gain a $\sqrt{\ep}$ in estimate (\ref{Erro1}).

By inspecting the proof of proposition \ref{Prop1}, we find out that the origin of the "bad" estimate is the term $\displaystyle \ep^{-1 }\pi^{j, \ep} \nabla \frac{\p P^D}{ \p x_j} $ in (\ref{2.4}). So we have to correct it in $\Oe_2$. Next,
the same type of difficulty arises with the  term $\displaystyle \ep^{-1}
 C^2_\pi  \frac{d}{d x_1} (\frac{\p P^D}{ \p x_2} |_\Sigma ) H(-x_2 )$ in   (\ref{2.14}). We handle it by modifying $\{ \mathbf{U}^\ep , P^\ep \}$. We include into the new velocity-pressure error pair the correction for the pressure term in (\ref{2.4}). The constant  $C^2_\pi$ corresponds to the behavior of $\omega^{2, bl}$ for $y_2 >0$ and we erase it in $\Oe_2$. Erasing it creates a pressure jump of order $O(\ep^{-1})$ and we compensate it by introducing a Darcy pressure field of such order in
$\Omega_2$.

We start by introducing an auxiliary problem correcting the singular pressure in (\ref{2.14}):
\begin{equation}\label{1.57newCell}
   \left\{ \begin{matrix}
   \hfill - \Delta_y \mathbf{w}^{i,k}_\pi (y) + \nabla_y \kappa^{i, k}_\pi (y ) = \pi^i (y)\mathbf{e}^k &
\hbox{ in } \ Y_F \cr \hfill \hbox{ div}_y \mathbf{w}^{i,k}_\pi (y) =0  & \hbox{ in
} \ Y_F \cr \hfill \mathbf{w}^{i,k}_\pi (y ) =0 & \hbox{ on } \, (\partial Y_F
\setminus
\partial Y) \cr
 \end{matrix}  \right.
\end{equation}
where $\pi^i$ is given by (\ref{1.57Cell}) and $\int_{Y_F} \kappa^{i, k}_\pi (y ) \ dy =0$.

Modified $\{ \mathbf{U}^\ep , P^\ep \}$ now read
\begin{gather}
\mathbf{{\tilde U}}^\ep = \mathbf{v}^\ep - H(x_2)  (\mathbf{u}^{eff} - \mathbf{e}^1  C^{2, bl}_1 \frac{\p P^D }{ \p x_2} |_\Sigma ) + H(-x_2)\{
 \sum_{j=1}^2 \mathbf{w}^{j, \ep} \frac{\p P^D }{ \p x_j}  - \notag \\
 \ep C^2_\pi \mathbf{w}^{1, \ep } \frac{d}{d x_1} \frac{\p P^D }{ \p x_2} |_\Sigma - \ep \sum_{j,k=1}^2 \mathbf{w}_\pi^{j,k} (\frac{x}{ \ep}) \frac{\p^2 P^D }{ \p x_j \p x_k} \}
-  \beta^{2,bl,\ep}  \frac{\p P^D }{ \p x_2} |_\Sigma ;\label{2.8A}\\
P^\ep =   p^{\ep } - H(x_2) p^{eff}   - H(-x_2) \{ \ep^{-2} P^D
-\,  \ep^{-1} ( C^2_\pi \frac{\p P^D }{ \p x_2} |_\Sigma + \sum_{j=1}^2 \pi^{j, \ep} \frac{\p P^D}{ \p x_j} )  +\notag \\
 \hskip-10pt C^2_\pi \pi^{1, \ep } \frac{d}{d x_1} \frac{\p P^D }{ \p x_2} |_\Sigma +\sum_{j,k=1}^2 \kappa_\pi^{j,k} (\frac{x}{ \ep}) \frac{\p^2 P^D }{ \p x_j \p x_k} )\} -\ep^{-1}  ( \pi^{2,bl,\ep} -C^2_\pi H(x_2 ) ) \frac{\p P^D }{ \p x_2} |_\Sigma .
 \label{2.9A}\end{gather}
Now all force-type terms are estimated as $\displaystyle C \sqrt{\ep}  \Vert \varphi \Vert _{H^1 (\Oe)^2 }$.
Furthermore, all normal stress jumps are of order $O(1)$.

Continuity of traces fails, but we correct it on the same way as in the original construction in subsection \ref{secpart}. The correction is of the order $O(\ep^{3/2})$ in $L^2$ for the velocity and of order $O(\sqrt{\ep})$ in $L^2$ for the pressure
and does not contribute to the result. Next we correct the effects on the boundary $\{ x_2 =-H \}$ and the compressibility effects. They are all of the next order and do not contribute to result.

The calculations yield the following estimates
\begin{gather}
|| \nabla \mathbf{{\tilde U}}^{2, \ep} ||_{L^2 (\Oe )^4} \leq C\sqrt{\ep}, \label{Erro1A}
   \\
  || \mbox{ div } \mathbf{{\tilde U}}^{2, \ep} ||_{L^2 (\Oe )^4}  \leq C\ep \label{Erro2A}
   \\
   ||  \mathbf{{\tilde U}}^{2, \ep} ||_{L^2 (\Oe_2 )^2} \leq C\ep^{3/2} \label{Erro2AA}
   \\
||  {\tilde P}^{2,\ep} ||_{L^2 (\Oe)} \leq \frac{C}{\sqrt{\ep}} \label{Erro3A}
\end{gather}
Now (\ref{Erro1A})-(\ref{Erro2AA}) imply
\begin{equation}\label{SigmaerrA}
  ||  \mathbf{{\tilde U}}^{2, \ep} ||_{L^2 (\Sigma )^2} \leq C\ep
\end{equation}
The rest of the proof is identical to the proof of theorem \ref{Thbasic}. Only difference is that we have gained a $\sqrt{\ep}$ in the estimates.\eop

\end{document}